\pdfoutput=1
\documentclass[review]{siamart171218}
\usepackage{titlesec}
\pdfoutput=1
\usepackage{soul}
\usepackage{enumitem}
\usepackage{lipsum}
\usepackage{amsfonts}
\usepackage{epstopdf}
\usepackage{multirow}
\usepackage[colorinlistoftodos,prependcaption,textsize=tiny]{todonotes}
\ifpdf
  \DeclareGraphicsExtensions{.eps,.pdf,.png,.jpg}
\else
  \DeclareGraphicsExtensions{.eps}
\fi

\usepackage{datetime}
\newdateformat{monthyeardate}{%
  \monthname[\THEMONTH] \THEDAY, \THEYEAR}

\makeatletter
\newcommand*{\addFileDependency}[1]{
  \typeout{(#1)}
  \@addtofilelist{#1}
  \IfFileExists{#1}{}{\typeout{No file #1.}}
}
\makeatother

\usepackage{academicons}
\usepackage{xcolor}

\usepackage{amsmath}
\allowdisplaybreaks
\usepackage{amssymb}
\usepackage{commath}
\usepackage{mathtools}
\usepackage{bbm}

\usepackage{color}
\usepackage{graphicx}
\usepackage[small]{caption}
\usepackage{subcaption}

\usepackage{relsize}
\usepackage{adjustbox}
\usepackage{algorithm,algpseudocode}
\usepackage{booktabs}
\usepackage{tikz}

\usepackage{verbatim}

\usepackage{enumitem}
\setlist[enumerate]{leftmargin=.5in}
\setlist[itemize]{leftmargin=.5in}

\newsiamthm{problem}{Problem}
\newsiamremark{remark}{Remark}
\newsiamremark{hypothesis}{Hypothesis} 
\crefname{hypothesis}{Hypothesis}{Hypotheses}
\newsiamremark{example}{Example}
\newsiamthm{claim}{Claim}
\newsiamthm{conjecture}{Conjecture}

\headers{Dynamic inverse problems with optical flow}{
Okunola T., Pasha M., Kilmer M., Freitag M.
}

\title{Efficient dynamic image reconstruction with motion estimation
\thanks{
\monthyeardate\today 
}}

\author{
Toluwani Okunola\thanks{Department of Mathematics, Tufts University} (\email{toluwani.okunola@tufts.edu})
\and 
Mirjeta Pasha\thanks{Department of Mathematics, Virginia Tech} (\email{mpasha@vt.edu})
\and
Misha E. Kilmer\thanks{Department of Mathematics, Tufts University}    (\email{misha.kilmer@tufts.edu})
\and
Melina Freitag\thanks{Institute of Mathematics, University of Potsdam} (\email{melina.freitag@uni-potsdam.de})
}

\usepackage{amsopn}

\usepackage{ifthen}
\usepackage{float}
\floatstyle{ruled}
\newfloat{algorithm}{htb}{alg}%[section]
\floatname{algorithm}{Algorithm}
\newcounter{algo@row}
\newcounter{algo@rowindent}
\newcommand{\algofont}[1]{\textbf{#1}}% S1
\newcommand{\algonumbersize}[1]{\scriptsize{#1}}% S2
\newcommand{\algopreitem}[1][\arabic{algo@row}]{\texttt{\algonumbersize{#1}}}
\newcommand{\algoitemskip}{\hspace{\value{algo@rowindent}cc}}
\newcommand{\algonewnestedopen}[2]{
	\newcommand{#1}[1][]{%
		\ifthenelse{\equal{##1}{}}{\item}{\item[{\algopreitem[##1]}]}
		\algoitemskip\algofont{#2}%
		\addtocounter{algo@rowindent}{1}%
		\ignorespaces
	}
}
\newcommand{\algonewnestedaux}[2]{
	\newcommand{#1}[1][]{
		\addtocounter{algo@rowindent}{-1}
		\ifthenelse{\equal{##1}{}}{\item}{\item[{\algopreitem[##1]}]}
		\algoitemskip\algofont{#2}%
		\addtocounter{algo@rowindent}{+1}%
		\ignorespaces
	}
}
\newcommand{\algonewnestedclose}[2]{
	\newcommand{#1}[1][]{
		\addtocounter{algo@rowindent}{-1}
		\ifthenelse{\equal{##1}{}}{\item}{\item[{\algopreitem[##1]}]}
		\algoitemskip\algofont{#2}%
		\ignorespaces
	}
}
\newcommand{\algonewcommand}[2]{
	\newcommand{#1}[1][default]{
		\ifthenelse{\equal{##1}{default}}{\item}{\item[{\algopreitem[##1]}]}%
		\algoitemskip\algofont{#2}%
		\ignorespaces
	}%
}
\newcommand{\algonewkeyword}[2]{\newcommand{#1}{\algofont{#2}}}
\algonewcommand{\STATE}{\ignorespaces}
\algonewcommand{\INPUT}{Input: }
\algonewcommand{\pINPUT}{\phantom{Input: }}
\algonewcommand{\COMPUTE}{Compute: }
\algonewcommand{\OUTPUT}{Output: }
\algonewcommand{\pOUTPUT}{\phantom{Output: }}

\algonewnestedopen{\IF}{if }
\algonewnestedaux{\ELSEIF}{else if }
\algonewnestedaux{\ELSE}{else }
\algonewnestedclose{\ENDIF}{end if }
\algonewnestedopen{\FOR}{for }
\algonewnestedclose{\ENDFOR}{end for }
\algonewnestedopen{\WHILE}{while }
\algonewnestedclose{\ENDWHILE}{end while }
\algonewcommand{\BREAK}{break}%

\algonewkeyword{\To}{to }%
\algonewkeyword{\Do}{do }%
\algonewkeyword{\Then}{then }%
\algonewkeyword{\End}{end }%
\algonewkeyword{\AND}{and }%
\algonewkeyword{\True}{true }%
\algonewkeyword{\False}{false }%
\algonewkeyword{\irbleigs}{irbleigs }%
\algonewkeyword{\tridiag}{tridiag}%
\algonewkeyword{\reorth}{reorth}%
\usepackage{cleveref}
\crefname{section}{§}{§§} % Singular and plural forms

\DeclareMathOperator*{\argmin}{arg\,min}

\newcommand{\R}{\mathbb{R}}

\newcommand{\bA}{{\bf A}}
\newcommand{\bB}{{\bf B}}

\newcommand{\bH}{{\bf H}}
\newcommand{\bI}{{\bf I}}

\newcommand{\bL}{{\bf L}}
\newcommand{\bM}{{\bf M}}

\newcommand{\bP}{{\bf P}}
\newcommand{\bQ}{{\bf Q}}
\newcommand{\bR}{{\bf R}}

\newcommand{\bU}{{\bf U}}
\newcommand{\bV}{{\bf V}}

\newcommand{\bb}{{\bf b}}

\newcommand{\br}{{\bf r}}
\newcommand{\bs}{{\bf s}}

\newcommand{\bu}{{\bf u}}
\newcommand{\bv}{{\bf v}}
\newcommand{\bw}{{\bf w}}

\newcommand{\by}{{\bf y}}
\newcommand{\bz}{{\bf z}}

\usepackage{textcomp}
\usepackage{gensymb}

\newcommand{\T}[1]{\boldsymbol{\mathcal{#1}}}

\newcommand{\bUpsilon}{{\boldsymbol{\Upsilon}}}

\usepackage{amsmath}
\usepackage{amssymb}
\usepackage[utf8]{inputenc}

\graphicspath{{figures_/}}
\begin{document}
\nolinenumbers
\maketitle

\begin{abstract}
Dynamic inverse problems are challenging to solve due to the need to identify and incorporate appropriate regularization in both space and time.  Moreover, the very large scale nature of such problems in practice presents an enormous computational challenge.    

In this work, in addition to the use of edge-enhancing regularization of spatial features, we propose a new regularization method that incorporates a temporal model that estimates the motion of objects in time.  In particular, we consider the optical flow model that simultaneously estimates the motion and provides an approximation for the desired image, and we incorporate this information into the cost functional as an additional form of temporal regularization.  We propose a computationally efficient algorithm to solve the jointly regularized problem that leverages a generalized Krylov subspace method.    
We illustrate the effectiveness of the prescribed approach on a wide range of numerical experiments, including limited angle and single-shot computerized tomography.
\end{abstract}

% REQUIRED


% Once the paper is published
\begin{DOI}
    Not yet assigned
\end{DOI}

\maketitle

\section{Introduction}
Classical X-ray tomography is a well-known technique used to explore and recover the unknown interior structure of a non-transparent static object of interest through the collection of data that are a result of the total attenuation of X-rays from different angles. The limited nature of the projection data and the presence of noise in the measured data means that regularization is required to ensure a unique solution that is less sensitive to noise \cite{EHN96, H98} .

In order to regularize an ill-posed problem suitable prior information is incorporated into the reconstruction. 
 
One way this can be accomplished is to include constraint(s) in the cost functional.  If $\bu$ denotes the vector of unknowns, then the regularized inverse problem results in the minimization problem
\[ \min_{\bu \in \R^{n}}{\mathcal{J}(\bu)} = \min_{\bu \in \R^{n}} \mathcal{D}(\bu) + \sum_{i} \lambda_{i} \mathcal{R}_i(\bu), \]
where $\mathcal{D}(\bu)$ denotes the data misfit term and the remainder denotes the regularization terms. 
 
For static image reconstruction, there is typically only one regularization term and parameter.  Common choices of regularization terms include the $\ell_p$ norm applied to $\bu$ or its gradient. Note that the regularization terms are premultiplied by  regularization parameters $ \lambda_{i}$.  The value of the parameter establishes the balance between the data fidelity term and the regularization term. Finding good values is a non-trivial problem and can pose a significant computational bottleneck in the reconstruction. See, for instance, \cite{ buccini2021linearized, pasha2023computational, pasha2023recycling, pashapythonpackage, chung2022, gazzola2020inner}. 

In our applications, however, the target of interest is dynamic. Dynamic image reconstruction (a.k.a. spatio-temporal reconstruction) refers to the process of recovering a sequence of (static) images that represent a time-varying scene. Thus, one needs to recover $n_t$ images simultaneously.
Moreover, the measured data are typically noisy and underdetermined in {\it both} space and time.
This time variation can lead to additional reconstruction artifacts if not corrected through the use of appropriate regularization \cite{mutaf2007impact, milanfar1999model}.
Thus, regularization in both space and time is needed, resulting in a very large scale inverse problem and the need to select more than one regularization parameter.  

With the ability of solving ever larger problems computationally, dynamic inverse problems have become increasingly popular \cite{Schuster_2018}, often using frameworks from data assimilation \cite{Stuart_2010}. Examples of dynamic inverse problems in tomography include imaging of moving objects or body parts, such as the heart or lungs \cite{nolte2022inverse}, dynamical impedance tomography \cite{schmitt1, schmitt2}, passive seismic tomography \cite{westman2012, Zhang2009PassiveST}, and undersampled dynamic x-ray tomography \cite{burger2017variational}. 
 
 To solve for all the frames at once, one needs to respect both spatial and temporal variations in the image sequence.

 For example, changes over time can be controlled through direct regularization in time (see  \cite{pasha2023computational, papoutsellis2021core, chung2018efficient}); through the use of low-rank methods to separate the dynamic particles on the object from the static background \cite{gao2011robust,arridge2022joint}; via Kalman filters and smoothing methods \cite{hakkarainen2019undersampled, bardsley2013krylov}, and more recently Bayesian methods \cite{chung2017generalized, lan2023spatiotemporal}.

The prevalence of \textit{dynamic inverse problems} argues for the development of models that describe the time evolution of images \cite{hakkarainen2019undersampled}. 

This leads to a more sophisticated class of regularization methods which will employ these models and use the corresponding estimates of frame-to-frame motion they provide to inform the temporal regularization. Examples of such reconstruction approaches include  diffeomorphic methods, which employ mathematical frameworks for smooth transformations \cite{chen2019new}, shape-based methods \cite{niemi2015dynamic, maboudi2024inferring, ozsar2025parametric}, and optical flow-based methods that utilize brightness constancy assumptions \cite{Burger2015OnOF, burger2017variational, schonelieb2018}. Wang et al. \cite{wang2023comprehensive} and Hauptmann et al. \cite{hauptmann2021image} provide extensive reviews of methods incorporating explicit motion models. 
Diffeomorphic methods are essential in applications where preserving the anatomical topology of structures is vital. 
Kadu et al. in \cite{kadu2024single} proposed a dynamic shape-sensing framework that employs level-sets in conjunction with spatio-temporal compression. Their method tracks the motion of the object by evolving its boundary and updating it as new image data becomes available. Although computationally efficient, their approach assumes that the movement of objects can be accurately represented by a small set of predefined functions, which may not always capture the complexity of real-world dynamics.
In the applications of interest to us in this paper, we assume no deformation of the object(s) of interest, so we will not consider diffeomorphic or shape-based methods further. 

We discretize in two space dimensions and consider the following \emph{state-space} model for the unknown state vector $\bu(t)$: 
\begin{align}
%\label{eq: generalProblem}
    \bu(t+1) & = \tilde{\mathcal{M}}\left(\bu(t)\right) + \xi^{t}, \quad \xi^{t} \sim \mathcal{N}\left(0, \boldsymbol{\Gamma}_t\right) \label{eq: dynamicEq}\\
    \bb^{t+1}  & = \bH^{t+1}\bu(t+1) + \gamma^{t}, \quad, \gamma^{t} \sim \mathcal{N}\left(0, \boldsymbol{\Sigma}_t\right), \quad t = 1, 2, \cdots, n_t. \label{eq: measurementEq} .
\end{align}
Here, $\bH^{t+1} \in \R^{m_s \times n_s}$ with $n_s = n_x \times n_y$ ($n_x$ and $n_y$ is the number of pixels on the x and y direction) is the observation operator at time $t$ that is also known as the \emph{parameter-to-observable} map whose application to the object of interest $\bu(t+1) \in \R^{n_x \cdot n_y}$ generates the observations $\bb^{t+1} \in \R^{m_s}$ with $m_s$ being defined as the products between the total number of angles and rays shined through the object. In the applications of interest to us, $m_s < n_s$ as a result of the limited angle computerized tomography. The vector $\gamma \in \R^{m_s}$ corresponds to Gaussian noise present in the measured data with covariance matrix $\boldsymbol{\Sigma}_t$. The model $\tilde{\mathcal{M}}$ represents dynamics that describe (along with the noise $\xi_t$) how images evolve from one time step to the next
 
In this paper, we focus on a specific type of state model, $\tilde{\mathcal{M}}$, that is known as optical flow.  Optical flow is a foundational concept in computer vision, describing the apparent motion of objects or textures between consecutive frames of an image sequence.  We describe this in more detail in Section \ref{subsec:optical_flow}. 

In the context of joint image reconstruction and motion estimation, optical flow plays a crucial role in estimating the motion field between consecutive frames. We can then use this information to reconstruct images from incomplete or degraded data. Several studies have proposed variational models that incorporate optical flow constraints to regularize the reconstruction process \cite{Burger2015OnOF, zhao_compensated, lucka2018enhancing, lucka2019palm}. These models typically alternate between estimating the optical flow, and solving the image reconstruction problem with a regularization term that promotes adherence of the reconstructed images to the estimated motion fields.

Despite the progress in spatio-temporal image reconstruction, several challenges remain, particularly in terms of computational efficiency. Methods that incorporate both motion estimation and image reconstruction are often computationally expensive, as they involve solving large-scale optimization problems with regularization terms. Methods based on proximal algorithms, such as the Alternating Direction Method of Multipliers (ADMM), as well as Primal Dual Methods have been successfully applied to solve regularized image reconstruction problems \cite{burger2017variational, Chambolle_Pock_2016, burger2014firstorderalgorithms}. These methods decompose the problem into simpler subproblems, which can be solved iteratively and in parallel, making them more computationally feasible. 

However, a major drawback of the existing methods (e.g. based on ADMM) is the sensitivity to several parameters, including (but not limited to) the regularization parameters. In applications, these parameters are often chosen by trial and error or on the basis of prior knowledge. Having to estimate good parameters adds a huge amount of computational overhead. And, as we demonstrate in the appendix of this paper, these methods are liable to exhibit sensitivity to parameters.  

In contrast, in MMGKS (Majorization Minimization Generalized Krylov Subspace)   
%package in Trips-Py 
when applied to solve a Tikhonov-type regularized problem, computes the regularized parameter cheaply and `automatically'.   

We propose an algorithm that alternates between recovery of the optical flow and the simultaneous reconstruction of all images in the sequence using the current estimate of the optical flow as a regularizer.  Both of these subproblems are regularization problems of the Tikhonov type, and therefore, we employ MMGKS to tackle each subproblem efficiently. 

Computationally, our algorithm is on par with other existing joint motion estimation methods and algorithms, but qualitatively, our method produces superior reconstructions. 

\paragraph{Overview of Main Contributions}
Our primary contributions are as follows:

\begin{itemize} \item We give a new, computationally efficient, spatio-temporal regularization algorithm, MMGKS-OF, to solve dynamic inverse problems in the setting when the objects move in time but are not deformed.    
Innovations include:  \begin{itemize} \item Solving the optical flow problem to estimate the motion model, \item Approximating the estimated motion model as a linear operator, 
\item Choosing regularization parameters cheaply and automatically, with no hand tuning
and \item Integrating the temporal dimension through solving the joint image reconstruction and motion estimation problem. \end{itemize} 
\item We conduct comprehensive numerical comparisons with competing methods, demonstrating the enhanced quality of reconstructed images with MMGKS-OF. \end{itemize}

\paragraph{Outline} This paper is structured as follows: \cref{sec:bg} provides a setup and preliminaries on the dynamic inverse problems setting and the optical flow. In \cref{sec:of_matrix} we detail our matrix-based encoding of the optical flow formulation and \cref{sec:mm-gks-of} introduces and explains the proposed MMGKS-OF algorithm. Numerical experiments are presented in  \cref{sec:experiments}, and concluding remarks and future directions are discussed in \cref{sec:conclusion}.
\section{Setup and Preliminaries} 
\label{sec:bg}

\subsection{Dynamic inverse problems}\label{ssec:background1}

Let us consider a sequence of two-dimension-al (2D) images, each represented by a matrix $\bU{(t)} \in \mathbb{R}^{n_x \times n_y}$, where $n_x$ and $n_y$ denote the number of rows and columns, respectively, and $t = 1,2,\dots, n_t$. We use the vectorization operator, denoted by vec, which stacks the columns of a matrix into a vector. Applying this operator to each image matrix $\bU{(t)}$, we obtain the vectorized image $\bu{(t)} = \text{vec}(\bU{(t)}) \in \mathbb{R}^{n_s}$, where $n_s=n_xn_y$is the total number of pixels in each image. We then concatenate the vectorized images $\bu{(t)}$ and apply the vectorization operator again: $\bu = \text{vec}([\bu{(1)}, \ldots, \bu{(n_t)}]) \in \mathbb{R}^{n}$, where $n=n_sn_t$ is the total number of pixels across all time frames, where $\bu$ provides a representation of the entire image sequence.

We are interested in solving inverse problems in space and time with an unknown target of interest. The goal is to recover from available measurements $\bb^{t} \in \mathbb{R}^{m_t}$, for $t = 1, 2,\dots, n_t$,
the {images} $\bu{(t)} \in \mathbb{R}^{n_s}$, {whose entries} represent pixels in the image. Since we focus on imaging applications, we use the term `pixels' (rather than `parameters') throughout the paper. 
Given the number of time points $n_t$, 
$m = \sum_{t=1}^{n_t} m_t$ is the total number of available measurements. We consider the number of pixels, $n_s$, to be fixed for all time points.

In dynamic inverse problems, the forward operator \(\mathbf{H} \in \mathbb{R}^{m \times n}\) may either evolve in time or remain constant. If time-dependent, \(\mathbf{H}\) is a block-diagonal matrix,
\[
\mathbf{H} = \text{blockdiag}(\mathbf{H}^{1}, \mathbf{H}^{2}, \dots, \mathbf{H}^{n_t}),
\]
where each block \(\mathbf{H}^{t} \in \mathbb{R}^{m_t \times n_s}\) represents the forward operator at time \(t\). If the forward operator is identical at all time steps, this leads to $\mathbf{H} = \mathbf{I}_{n_t} \otimes {\mathbf{H}}^{t},$ where \(\otimes\) denotes the Kronecker product. The vector $\bb ={\rm{vec}}([\bb^{1},\dots,\bb^{n_t}])\in \R^{m}$ represents measured data that are contaminated by an unknown error {(or noise)} $\gamma \in \R^{m}$ 
that may stem from measurement errors. 

The end goal of the inverse problem involves recovering the pixel values of all images over all time (i.e. $\bu$) from the data $\bb$ along with the parameters of the state model $\mathcal{M}$, where $\mathcal{M}(\bu)^T = [\tilde{\mathcal{M}}(\bu(1)),\tilde{\mathcal{M}}(\bu(2)),\ldots,\tilde{\mathcal{M}}(\bu(n_t))]^T$. We formulate this as a minimization problem
\begin{equation}\label{eq: l2-lq-general}
\min_{\bu \in \R^{n}}{\mathcal{J}(\bu)} := \mathcal{D}(\bu)+\lambda_1\mathcal{R}(\bu)+\lambda_{2} \mathcal{M}(\bu),
\end{equation}
where the functional $\mathcal{D}(\bu)$ that represents the data-fidelity term that, if the measurement noise is assumed to be Gaussian, takes the form $\frac{1}{2}\|\bH\bu-\bb\|_2^2$. {Parameters $\lambda_j > 0$, $j=1,2$ are known as regularization parameters that determine the balance between the data-misfit and the regularization terms $\mathcal{R}(\bu)$ and $\mathcal{M}(\bu)$.} The functional $\mathcal{R}(\bu)$ is a regularization term that can take different forms. One of the most popular and general forms is $\|\mathbf{\Psi} \bu\|_q$ for some choice of the regularization matrix $\boldsymbol{\Psi} \in \R^{r \times n}$ and $0 < q \leq 2$. Frequently, a selection of $q = 1$ is made to ensure adequate representation of edges and that the regularization term defines a valid norm \cite{buccini2021linearized, chung2019flexible, chung2022, lanza2015generalized, buccini2020modulus}. We therefore use $q = 1$ throughout the paper, but we note that {\it the machinery we develop here can easily be adapted to other choices of $q$}. 
The functional $\mathcal{M}(\bu)$ encodes temporal knowledge through a model $\mathcal{M}$. Details on the choice of the model and techniques to estimate its parameters are presented in the following sections.

\subsection{Optical Flow}\label{subsec:optical_flow}

Optical flow is defined as the apparent motion between a sequence of images. As early as 1981, a method for computing optical flow was introduced in \cite{HORN1981185}. Since then, it has been widely used as a technique for describing image motion, particularly for image sequences with small time differences between consecutive frames and has been significantly useful in computer vision \cite{alfarano2024estimating}.

An optical flow rests on the assumption that the so called \emph{optical flow constraint} (OFC) is satisfied. This constraint enforces that as the pixels in the image move or change position with time, their intensities remain the same. In particular, given an image $\bu$ and considering a single pixel at location $(x^i,y^i)$, the optical flow constraint states 
% \textcolor{red}{Write * with i}.
\begin{equation}
\frac{d\mathrm{\bu} (x^i,y^i,t)}{dt} = 0.
\end{equation}

We are interested in modeling the \emph{state model} $\mathcal{M}$ as optical flow, so we assume that the image sequences satisfy the optical flow assumption. Hence, we consider the following dynamic for the motion 
 of the pixels in the images: considering a single pixel $i$ in $\bu(t)$ located at $(x^i(t),y^i(t))$ with intensity $\bu(x^i(t),y^i(t))$. If this pixel has velocity $\bs^i(t) := (s_x^i(t),s_y^i(t))$, then at time $t+ \Delta t$, the pixel location changes to $(x^i(t) + s_x^i(t) \Delta t, y^i(t) + s_y^i(t) \Delta t)$ (we assume that the velocity is constant in the interval $[t,t+\Delta t)$), but its intensity stays the same. Hence, the following holds
 \begin{equation}\label{eq:opt_flow_2}
\bu(x^i + s_x^i(t) \Delta t, y^i + s_y^i(t) \Delta t, t+ \Delta t) = \bu(x^i, y^i,t).
 \end{equation}
Since $\Delta t$ is small, we linearly approximate the RHS of \eqref{eq:opt_flow_2} around $(x^i, y^i,t)$
\begin{align}\label{eq: linearization}
\bu(x^i(t) + s_x^i(t) \Delta t, y^i(t) + s_y^i(t) \Delta t, t+ \Delta t) \approx 
\bu(x^i(t), y^i(t),t) + &\\ \Delta t \left(s_x^i(t) \frac{\partial \bu}{\partial x} (x^i(t), y^i(t),t) +  s_y^i(t) \frac{\partial \bu}{\partial y}(x^i(t), y^i(t),t) + \frac{\partial \bu}{\partial t} (x^i(t), y^i(t),t)\right).\nonumber
\end{align}
Substituting \eqref{eq: linearization} 
%for  $\bu(x^i + s_x^i(t) \Delta t, y^i + s_y^i(t) \Delta t, t+ \Delta t)$ 
in \eqref{eq:opt_flow_2} gives 
\begin{align*}
\bu(x^i(t), y^i(t),t) + \\
\Delta t \left( s_x^i(t) \frac{\partial \bu}{\partial x}(x^i(t), y^i(t),t) + s_y^i(t) \frac{\partial \bu}{\partial y}(x^i(t), y^i(t),t) +\frac{\partial \bu}{\partial t}(x^i(t), y^i(t),t)\right) \\ = \bu(x^i(t), y^i(t),t),
\end{align*}
equipping us with the equation
\begin{equation}\label{eq: of}
 \frac{\partial \bu}{\partial x}(x^i(t), y^i(t),t) s_x^i(t) + \frac{\partial \bu}{\partial y}(x^i(t), y^i(t),t) s_y^i(t)  + \frac{\partial \bu}{\partial t}(x^i, y^i,t) = 0,
\end{equation}
For the rest of the paper, to simplify notation, we agree to denote 
\[\bu_x^i(t)= \frac{\partial \bu}{\partial x}(x^i(t), y^i(t),t), \quad \bu_y^i(t)= \frac{\partial \bu}{\partial y}(x^i(t), y^i(t),t),\quad \bu_t^i(t)= \frac{\partial \bu}{\partial t}(x^i(t), y^i(t),t),\]
which leads to the following equation
% Finally, we can write \eqref{eq: of} as  
\begin{equation*}
\begin{pmatrix}
    \bu_x^i(t) & \bu_y^i(t)
\end{pmatrix}
\begin{pmatrix}
    s_x^i(t) \\ s_y^i(t)
\end{pmatrix}
= - \bu_t^i(t).
\end{equation*}

For each pixel $i = 1, 2, \dots, n_s$ in the image $\bu(t)$, we can compute the velocities of $\bu$ at time $t$, and obtain the system of equations:
\begin{equation*}
\underbrace{\begin{pmatrix}
    \bu_x^i(t) & \ldots & 0 & \bu_y^i(t) & \ldots & 0 \\
    & \ddots & & & \ddots & \\
  0 & \ldots &\bu_x^{n_s}(t)&0 & \ldots & \bu_y^{n_s}(t)
\end{pmatrix}}_{\bUpsilon(\bu(t))}
\underbrace{\begin{pmatrix}
    s_x^1(t)\\ \vdots \\ s_x^{n_s}(t)\\ s_y^1(t)\\ \vdots    \\ s_y^{n_s}(t)
\end{pmatrix}}_{\bs(t)} = - 
\underbrace{\begin{pmatrix}
\bu_t^i(t) \\ \vdots \\ \\ \vdots \\
\bu_t^{n_s}(t)
\end{pmatrix}}_{\bu_t}
\end{equation*}

This is an inverse problem and its solution can be obtained by solving the least-squared problem
\begin{equation}
    \min _ {\mathbf{s} \in R^{2n_s}} \| \bUpsilon(\bu(t))\mathbf{s}(t) + \bu_t(t)\|_p^p.
    \label{eq:of_lsq}
\end{equation}

However this equation is under-determined and could have many solutions, so an additional constraint is usually imposed on the velocity vectors \cite{HORN1981185} to produce the regularized problem:
\begin{equation}
    \min _ {\mathbf{s} \in R^{2n_s}} \| \bUpsilon(\bu(t))\mathbf{s}(t) + \bu_t(t)\|_p^p + \gamma \left\| \begin{pmatrix}\nabla 
        \mathbf{s}_x(t) \\ 
        \nabla \mathbf{s}_y(t)
    \end{pmatrix}\right\|_q^q,
    \label{eq:of_reg}
\end{equation}

where $\mathbf{s}_x(t) = \begin{pmatrix} s_x^1(t) & \cdots & s_x^{n_s}(t) \end{pmatrix}^\top$, $\mathbf{s}_y(t) = \begin{pmatrix} s_y^1(t) & \cdots & s_y^{n_s}(t) \end{pmatrix}^\top$,  $\gamma \in \mathbb{R}$ is and a regularization parameter, and  $p$ and $q$ can take values in $\{1,2\}$.

For $\bs = \{\mathbf{s}(j)\}_{j=1}^{n_t}$, we can write the problem jointly as:
\begin{equation}
    \min _ {\mathbf{s} \in R^{2n_s}} \| \bUpsilon(\bu) \mathbf{s} + \bu_t\|_p^p  + \gamma \|
        \hat{\bL} \mathbf{s}\|_q^q,
    \label{eq:of_reg_full}
\end{equation}
where $\bUpsilon(\bu) = \text{diag}(\bUpsilon(\bu(1)), \ldots, \bUpsilon(\bu(n_t - 1)))$, $\bu_t = \begin{pmatrix} \bu_t(1) & \cdots & \bu_t(n_t-1)\end{pmatrix}^\top$, and $\hat{\bL}$ is a suitable regularization operator representing the discrete gradients.

A choice of $p=2$ can be interpreted as a constraint imposing smoothness on the velocity vector, while choosing $q=1$ imposes sparsity/ piece-wise constancy of the velocity vectors.

\section{Encoding Optical Flow Information in a Matrix}
\label{sec:of_matrix}
To be able to use the optical flow information as regularization, we first need to discuss how to 
encode the information into a matrix. 

We define the motion function $\mathcal{M}: \mathbb{R}^{n_x \times n_y} \times \mathbb{R}^{n_x \times n_y\times 2} \to \mathbb{R}^{n_x \times n_y}$ for an image sequence. Given an image at time $t + \Delta t$, $\mathbf{u}(t + \Delta t)$, and the velocity at time $t$, $\mathbf{s}(t)$, $\mathcal{M}$ returns the image at time $t$, $\mathbf{u}(t)$, satisfying the optical flow constraint \eqref{eq:opt_flow_2}:
\[
\mathcal{M}(\mathbf{u},\mathbf{s}) := \tilde{\mathbf{u}},
\]
where $\tilde{\mathbf{u}}$ satisfies $\tilde{\mathbf{u}}(x,y,t) = \mathbf{u}(x+s_x,y+s_y,t+\Delta t)$. We refer to $\mathcal{M}$ as the \emph{reverse motion}, representing movement from time $t + \Delta t$ to $t$. The forward motion, from $t$ to $t + \Delta t$, uses velocity vectors $\mathbf{s}'(t+\Delta t)$. Incorporating this reverse motion is useful for solving dynamic inverse problems.

The motion model $\mathcal{M}$ is not necessarily linear in $\bu$. To incorporate it as a linear regularization term in the MMGKS algorithm, we express it as a linear function of $\mathbf{u}$. 

In an ideal setting, since the optical flow equation preserves pixel intensity, $\mathcal{M}$ can be seen as a permutation based on velocity $\mathbf{s}$. After linearization, we denote the linear operator as $\mathbf{M}$ (an $n_s \times n_s$ matrix), such that:
\begin{equation}
    \mathbf{M}(\mathbf{s})= \left(\delta_{r,r+\mathbf{s}^r}\right)_{r=(1,1)}^{(n_x, n_y)},
    \label{eq:m_matrix}
\end{equation}
where $\delta$ is the Kronecker delta, and
\[
\mathcal{M}(\mathbf{u},\mathbf{s}) = \mathbf{M}(\mathbf{s}) \mathbf{u}.
\]
However, this is not usually the case, as many times, solving the optical flow equations gives continuous velocities. Given $\mathbf{u}(t+1)$ and optical flow $\mathbf{s}(t)$, we seek $\mathbf{M}(\mathbf{s})$ such that $\mathbf{M}(\mathbf{s}(t)) \mathbf{u}(t+1) \approx \mathbf{u}(t)$. This maps each pixel at $r = (x,y)$ to $r + \mathbf{s}^r = (x+s_x,y+s_y)$.
We consider two approaches:

\begin{enumerate}
\item \text{Rounding and Clipping:} We round velocity components to the nearest integers and clip source coordinates to image boundaries. The permutation matrix $\mathbf{M}(\mathbf{s}) \in \mathbb{R}^{n_s \times n_s}$ is then defined by Equation~\eqref{eq:m_matrix}.

\item \text{Bilinear Warping:} We construct M(s) as a bilinear warping operator, similar to the method in \cite{burger2017variational}. Each row of M(s) represents a target pixel, with its non-zero entries encoding the corresponding bilinear interpolation weights from the original image's neighboring pixels.
\[
\mathbf{M}(\mathbf{s})_{i, j} =
\begin{cases}
w_{x'y'}(j), & \text{if } j \in \mathcal{N}(r + \mathbf{s}^r), \\
0, & \text{otherwise},
\end{cases}
\]
where $i = y n_x + x$ is the row index corresponding to the target pixel at location $r = (x, y)$. The set $\mathcal{N}(r + \mathbf{s}^r)$ denotes the indices of the four nearest neighbors of the displaced location $r' = r + \mathbf{s}(r) = (x', y')$, where $x' = x + s_x$ and $y' = y + s_y$, and $w_{x'y'}(j)$ are the bilinear interpolation weights associated with these neighbors.
The column index $j$ in $\mathbf{M}(\mathbf{s})_{i, j}$ corresponds to the index of one of these four neighboring pixels in the original image.
\end{enumerate}
Given the encoding matrix $\mathbf{M}(\mathbf{s}(t))$, it holds that for an image $\mathbf{u}(t), \mathbf{u}(t+1)$ satisfying the brightness constancy assumption, 
\begin{align*}
(\mathbf{M}(\mathbf{s}(t))\mathbf{u}(t+1) - \mathbf{u}(t))(x,y) &= \mathbf{u}(x+\mathbf{s}_x, y+\mathbf{s}_y, t+1) - \mathbf{u}(x, y, t) \\
&\approx \mathbf{u}_x(x, y, t) \mathbf{s}_x (t) + \mathbf{u}_y(x, y, t) \mathbf{s}_y(t) + \mathbf{u}_t(x, y, t)\\
&= (\mathbf{\Upsilon}(\bu(t))\mathbf{s}(t) + \mathbf{u}_t)(x,y).
\end{align*}
for each pixel located at $(x,y)$.
If we let
\begin{align}
\bar{\mathbf{M}}(\mathbf{s}) &:= \begin{pmatrix}
    I & - \mathbf{M}(\bs(1)) & \mathbf{0}& \cdots & \mathbf{0} & \mathbf{0} \\
    \mathbf{0} & I & - \mathbf{M}(\bs(2)) & \mathbf{0} & \cdots & \mathbf{0} \\
    \vdots & \vdots & \vdots & \vdots & \vdots & \vdots \\
    \mathbf{0} & \mathbf{0} & \cdots & \mathbf{0} & I & - \mathbf{M}(\bs(n_t-1))
\end{pmatrix}.
\end{align}
Then, 
\begin{equation}
\bar{\mathbf{M}}(\mathbf{s}) \mathbf{u} \qquad \text {and } \qquad \mathbf{\Upsilon}(\bu) \mathbf{s} + \mathbf{u}_t
\end{equation}
represent the same term -- linear in $\mathbf{u}$ and $\mathbf{s}$, respectively.\newline
For the reverse optical flow, we let
\begin{equation}
\bar{\mathbf{M}'}(\mathbf{s}') := \begin{pmatrix}
- \mathbf{M}(\bs'(2))& I & \mathbf{0} & \mathbf{0}& \cdots & \mathbf{0} \\
\mathbf{0} & - \mathbf{M}(\bs'(3)) & I &\mathbf{0} & \cdots & \mathbf{0} \\
\vdots & \vdots & \vdots & \vdots & \vdots & \vdots \\
\mathbf{0} & \mathbf{0} & \cdots & \mathbf{0} & - \mathbf{M}(\bs'(n_t)) & I \\
\end{pmatrix}.
\end{equation}
Similarly, we have the following equivalence
\begin{equation}
 \bar{\mathbf{M}'}(\mathbf{s}') \mathbf{u}, \qquad \equiv \qquad  \mathbf{\Upsilon}'(\bu)  \mathbf{s}' + \mathbf{u}_t',   
\end{equation}
where $\mathbf{\Upsilon}'(\bu) = \text{diag}(\mathbf{\Upsilon}(\bu(2)), \ldots, \mathbf{\Upsilon}(\bu(n_t)))$, $\mathbf{u}_t' = \begin{pmatrix} -\mathbf{u}_t(2) & \cdots & -\mathbf{u}_t(n_t)\end{pmatrix}^\top$. 

This formulation is used in the following section.
\begin{remark}[On Estimating the Reverse Optical Flow]
    Following the OFC, we see that we can estimate $\mathbf{s}'$ from $\mathbf{s}$ using the following relation:
\begin{equation*}
    ({s}'_{x + s_x(t)}(t + \Delta t), {s}'_{y + s_y(t)}(t + \Delta t)) = - (s_x(t),s_y(t)),
    \label{eq:s_prime_from_s}
\end{equation*}
where $(s_x(t),s_y(t))$ is the velocity of a pixel at location $(x,y)$.
This is useful for saving computational expenses of solving the optical flow equations in the reverse direction, and we employ this relation in this paper to estimate $\mathbf{s}'$.
\end{remark}

\section{Majorization Minimization Generalized Krylov Subspace with Optical Flow Estimation (MMGKS-OF)}
\label{sec:mm-gks-of}

As we linearize the state model $\mathcal{M}$, we obtain a regularization matrix $\hat{\bM}\left(\bs, \bs'\right)$ that we can include as a regularization operator. In particular, in addition to the data fidelity term and an $\ell_1$ regularization term (a), the regularization term that is concerned with the velocity parameters (b) is being considered in the functional $\mathcal{J}(\bu, \bs, \bs')$,
\begin{align}
    \mathcal{J}(\bu, \bs, \bs') := \underbrace{\|\mathbf{H}\bu - \mathbf{b}\|_2^2 +  \lambda_1 \|\boldsymbol{\Psi}\bu\|_1}_{(a)} + \underbrace{\lambda_{2}\|\hat {\bM}(\bs,\bs')\bu\|_2 + \gamma\| \hat{\bL} \bs \|_2^2 + \gamma' \| \hat{\bL} \bs' \|_2^2}_{(b)}
\label{eqjointobj}, 
\end{align} where, recall that $\hat{\bL}$ represents a discrete gradient operator and 
\[\hat{\mathbf{M}}(\bs,\bs'): = \begin{bmatrix}
    \bar{\bM}(\bs) & \bar{\bM^{'}}(\bs')
\end{bmatrix}^\top.\]
This formulation appears to show four different regularization parameters, which might seem problematic.  However, in practice we can set $\lambda = \lambda_{1}= \lambda_2$, and as the regularization terms end up only contributing to one of each of the subproblems, we will use MMGKS to select those parameters automatically. 

The approach that we follow here is that of a \emph{joint image reconstruction and motion estimation} problem that seeks to reconstruct both the sequence of the images and the parameters of the state model. This leads to the following minimization problem
\begin{align}\label{eqjointproblem}
    \bu^{*}, \mathbf{s}^{*}, \bs'^{*} = \operatorname*{argmin}_{ \substack{\bu \in \mathbb{R}^{n}, \mathbf{s}, \mathbf{s}^{'} \in \mathbb{R}^{2n_s}}}
    \mathcal{J}(\bu, \bs, \bs'),
\end{align}
To solve \eqref{eqjointproblem}, we 
split the problem into two by alternating between the image reconstruction and optical flow estimation.
\paragraph{Initialization:} 
We begin with an initial estimate $\bu^{(0)}$ of the image sequence then we proceed with the estimation of the parameters of the optical flow.
\paragraph{Optical Flow Estimation:}
Given $\bu^{(k)}$, we fix $\bUpsilon^{(k)} := \bUpsilon(\bu^{(k)})$ and obtain $\bs^{(k)}$  and $\bs'^{(k)}$ by solving the optical flow equations:
\begin{equation}%\label{eq: OF_forward}
   \mathbf{s}^{(k)} = \operatorname*{argmin}_{\mathbf{s} \in R^{2n}} \| \bUpsilon^{(k)} \mathbf{s} + \bu_t^{(k)}\|_p^p + \gamma \| 
        \hat{\bL} \mathbf{s}\|_q^q, 
\label{eq:subproblem_motion_estimation(a)}
\end{equation}
\begin{equation}
   \mathbf{s'}^{(k)} \operatorname*{argmin}_ {\mathbf{s'} \in R^{2n}} \| \bUpsilon{'}^{(k)} \mathbf{s'} + \bu_t'^{(k)}\|_p^p + \gamma ' \| 
        \hat{\bL} \mathbf{s'}\|_q^q, 
\label{eq:subproblem_motion_estimation(b)}
\end{equation}

\paragraph{Image Reconstruction:}
Given $\mathbf{s}^{(k)}$, $\mathbf{s'}^{(k)}$, fix $\hat{\bM}^{(k)}:= \hat{\bM} (\mathbf{s}^{(k)}, \mathbf{s'}^{(k)})$, and solve the system:
\begin{equation}
    \bu^{(k+1)} = \operatorname*{argmin}_{{\bu \in \mathbb{R}^{n}}}
    \|\mathbf{H}\bu - \mathbf{b}\|_2^2 +  \lambda_1 \|\boldsymbol{\Psi}\bu\|_1 + \lambda_{2} \|\hat {\bM}^{(k)}\bu\|_1.
\label{eq:rec_problem_s_fixed} 
\end{equation}

\subsection{Numerical Implementation and Computational Costs}
Here, we provide a detailed explanation of our approach to solving the Image Reconstruction and Optical Flow Estimation sub-problems introduced earlier. We also elaborate on the methods inherent to the MMGKS algorithm (carried over to the MMGKS-OF algorithm) for selecting the regularization parameter. Subsequently, we highlight the computational costs associated with our proposed algorithm and discuss strategies employed to minimize these costs.

We let $\bL_{d} \in \R^{(n_d-1) \times n_d}$
be a rescaled finite difference discretization of the first derivative operator. Operators of this kind are known to damp fast oscillatory components of a vector $\bu^{(t)}$; see, for instance, a discussion in \cite{donatelli2014square}. We then define 
\begin{equation}\label{eq: L}
    {\bL} = \begin{bmatrix}
    \bI_{n_x} \otimes \bL_{y} \\ \bL_{x} \otimes \bI_{n_y}
\end{bmatrix},
\end{equation} 
where the matrix $\bL_x$ represents the discrete 
first derivative operator in the $x$-direction (vertical) and
$\bL_y$ in the $y$-direction (horizontal).

\subsubsection{Solving the  Image Reconstruction Subproblem}
To solve problem \eqref{eq:rec_problem_s_fixed}, we set $\lambda = \lambda_1 = \lambda_{2}$ to reduce the number of computed parameters. We also choose the regularization operator $\boldsymbol{\Psi} = \text{blockdiag}(\bL \ldots \bL)$. We then combine the regularization term and the motion fidelity term to obtain the problem: 
\begin{equation}\label{eq:rec_problem_s_fixed_joint}
    \bu^{(k+1)} = \operatorname*{argmin}_{{\bu \in \mathbb{R}^{n}}}
    \|\mathbf{H}\bu - \mathbf{b}\|_2^2 +  \lambda \|\boldsymbol{\Theta}^{(k)}\bu\|_1,
\end{equation}
where 
\[\boldsymbol{\Theta}^{(k)} = \begin{bmatrix}
    \boldsymbol{\Psi} \\ \hat{\bM}^{(k)}
\end{bmatrix}.\]
To solve \eqref{eq:rec_problem_s_fixed_joint}  we employ the MMGKS method \cite{lanza2015generalized} which proceeds as follows. 
Let $\bu^{(k)}$ be an approximate solution of \eqref{eq:rec_problem_s_fixed_joint}, we define $\bz^{(k)} = \boldsymbol{\Theta}^{(k)}\bu^{(k)}$ and the weighting matrix $\bP_{\epsilon}^{(k)} = \left(\text{diag}\left( \frac{1}{\phi_\epsilon(\bz^{(k)})}\right)\right)^{1/2},$
% \end{equation} 
where operations on the right are performed element-wise. Then we obtain the quadratic tangent majorant $\mathcal{Q}(\bu, \bu^{(k)})$ for $\mathcal{J}_{\epsilon,\lambda}(\bu)$,  
with weighting matrix 
$\bP_{\epsilon}^{(k)}$ as 
%as
\begin{equation}\label{eq: QuadraticMajorantQ}
\begin{array}{rcl}
\mathcal{Q}(\bu, \bu^{(k)})  &:=&
\|\bH\bu-\bb\|^{2}_{2}
+\lambda \left(\|\bP_{\epsilon}^{(k)}\boldsymbol{\Theta}^{(k)}\bu\|^{2}_{2}\right)+c,
\end{array}
\end{equation}
where $c$ denotes a suitable constant\footnote{A non-negative value of $c$ is technically required for the second condition in the definition to hold. However, since the minimizer is independent of the value of $c$, we will not discuss it further.} that is independent of $\bu^{(k)}$.
Minimizing the quadratic tangent majorant \eqref{eq: QuadraticMajorantQ} yields the next iterate $\bu^{(k+1)}$.
We proceed by setting the gradient of 
\eqref{eq: QuadraticMajorantQ}
to zero which gives the normal equations
\begin{equation}\label{eq: normaleqQuadMajorant}
(\bH^T\bH + \lambda  {\boldsymbol{\Theta}^{(k)}}^{T} (\bP_{\epsilon}^{(k)})^2\boldsymbol{\Theta}^{(k)})\bu = \bH^T\bb.
\end{equation}
The condition \(\mathcal{N}(\bH^T\bH)\cap \mathcal{N}({\boldsymbol{\Theta}^{(k)}}^{T} (\bP_{\epsilon}^{(k)})^2\boldsymbol{\Theta}^{(k)})=\{0\}\), which guarantees that the system \eqref{eq: normaleqQuadMajorant} has a unique solution, typically holds in practice. In that case, for fixed $\lambda$,
the solution  
of \eqref{eq: normaleqQuadMajorant} is the unique minimizer
of the quadratic tangent majorant $\mathcal{Q}(\bu, \bu^{(k)})$.
In practice, $\lambda$ is usually not known in advance, and we discus its computation in \cref{sec:regparam}.

Unfortunately, solving \eqref{eq: normaleqQuadMajorant} for large $\bH$ and  $\Theta^{(k)}$ may be computationally prohibitive.  Therefore, 
in \cite{huang2017majorization} the authors propose the 
MMGKS method that computes approximations by projection onto low dimensional subspaces. 
The method starts with a few steps of Golub-Kahan
bidiagonalization (GKB) 
with $\bH$ and $\bu_1 = \bb/\| \bb \|_2$ to determine the initial subspace $\bV_\ell$ such that 
\(\bH\bV_{\ell}=\bU_{\ell+1}\bB_{\ell}\),
and $\bV_\ell$ and
$\bU_{\ell+1}$ have orthonormal columns. 
Given $\bu^{(k)}$ and 
$\bP_{\epsilon}^{(k)}$ with $k=0$ in the first step), 
MMGKS computes the thin QR factorizations 
\begin{align}\label{eq: QR}
\bH\bV_{\ell+k} = \bQ_{\bH}\bR_{\bH}, \quad
\bP_{\epsilon}^{(k)}\Theta^{(k)}\bV_{\ell+k} = \bQ_{\Theta^{(k)}}\bR_{\Theta^{(k)}} .
\end{align}
which are now inexpensive to compute because of the smaller dimension. Restricting 
$\bu^{(k+1)}$ to 
${\rm range}(\bV_{\ell})$, 
$\bu^{(k+1)} = \bV_{\ell} \bz^{(k+1)}$, 
and adjusting
\eqref{eq: normaleqQuadMajorant}
accordingly leads to the
following small system
of equations for $\bz^{(k+1)}$, 
\begin{equation} \label{eq: minKryov2}
  (\bR_\bH^T\bR_\bH + 
  \lambda \bR_{\Theta^{(k)}}^T\bR_{\Theta^{(k)}})\bz^{(k+1)} =
  \bR_\bH^T\bQ_\bH^T\bb .
\end{equation}

Its solution gives the residual vector of the (full) normal equations:
\begin{equation*}\label{eq: residual}
\br^{(k+1)}=\bH^T(\bH\bV_{\ell}\bz^{(k+1)} -\bb)+\lambda  {\boldsymbol{\Theta}^{(k)}}^T(\bP_{\epsilon}^{(k)})^2\boldsymbol{\Theta}^{(k)}\bV_{\ell}\bz^{(k+1)}.
\end{equation*}
We expand the solution subspace with a new normalized residual
$\bv_{\rm new}=\frac{\br^{(k+1)}}{\|\br^{(k+1)}\|_2}$,
\begin{equation}\label{eq: enlargeMM-GKS}
\bV_{\ell+1}=[\bV_{\ell},\bv_{\rm new}]\in\R^{n\times(\ell+1)}.\end{equation} 
In exact arithmetic,  $\bv_{\rm new}$ is orthogonal to the columns of $\bV_{\ell}$, but in computer arithmetic reorthogonalization of $\bV_{\ell+1}$ is typically needed.
Next, $\bP_{\epsilon}^{(k)}$ is updated for the new solution estimate, and the process in steps \eqref{eq: QR} - \eqref{eq: enlargeMM-GKS} is repeated, expanding the solution space until a sufficiently accurate solution is reached. The MMGKS algorithm is summarized in Algorithm \ref{alg:mm-gks}  in  Appendix \ref{sec:mmgks_alg} and more details on the majorization and minimization steps can be found in \cite{lanza2015generalized}.

\begin{algorithm}[ht]
\caption{MMGKS-OF Algorithm}
\label{alg: mm-gks-dyn}
\begin{algorithmic}[1]
\Require: $\bH$, $\bL$, $\bb$, ${\bu}^{(0)}$, $\epsilon$
\Ensure: An approximate solution $\bu^i$
\STATE \textbf{function} \textsc{MMGKS-OF}{($\bH$, $\bL$, $\bb$,$\bu^{(0)}$,$\epsilon$)}
\STATE Generate the initial subspace basis $\bV_\ell \in \R^{n \times \ell}$ such that $\bV_\ell^T\bV_\ell = \bI $
\FOR {$k = 0, 1, 2, \ldots$ until convergence}
\STATE  $\mathbf{s}^{(k)}$ and $\mathbf{s}'^{(k)} = \{\textsc{SOLVE-OF} ({\bu^{(k)}(t), \bu^{(k)}(t+1), \bL})\}_{t=1}^{n_t-1}$  \Comment{Solve the optical flow equations  in both directions to obtain the velocity vectors}
\STATE $\hat{\mathbf{M}}^{(k)} = \hat{\bM}(\mathbf{s}^{(k)}, \mathbf{s}'^{(k)})$
\STATE $\Theta^{(k)} = \begin{pmatrix}
    \boldsymbol{\Psi} \\ \hat{\mathbf{M}}^{(k)}
\end{pmatrix}$
\STATE $\nu^{(k)} =  \Theta^{(k)} \bu^{(k)}$
% \STATE $\sigma_\epsilon^{(k)} = ((\nu^{(k)})^2 + \epsilon^2)^{-1/2}$
\STATE $\bP_\epsilon^{(k)}= (diag((\nu^{(k)})^2 + \epsilon^2))^{-1/4}$
\STATE Obtain $\bH\bV_{\ell+k}$, and $\bP_\epsilon^{(k)}\Theta \bV_{\ell+k}$
\STATE Compute their QR decompositions: $\bH\bV_{\ell+k} = \bQ_{\bH}\bR_{\bH}$, $\bR_\epsilon^{(k)}\Theta \bV_{\ell+k} = \bQ_\Theta \bR_\Theta$
\STATE Optimally select $\lambda^{(k)}$.
\STATE Solve the least squares problem \eqref{eq: minKryov2} with $\lambda^{(k)}$ to obtain $\bz^{(k+1)}$
\STATE $\bu^{(k+1)} = \bV_{\ell+k}\bz^{(k+1)}$ 
\STATE $\br^{(k+1)} = \bH^T(\bH\bV_{\ell+k}\bz^{(k+1)} - \bb) +\lambda^{(k)}  \mathbf{\Theta}^{{{(k)}^T}}\bP_\epsilon^{(k)}\mathbf{\Theta}^{(k)}\bV_{\ell+k}\bz^{(k+1)} $ 
\STATE $\br^{(k+1)}= \br^{(k+1)} - \bV_{\ell+k}^T \br^{(k+1)}$
\STATE $\bv = \tfrac{\br^{(k+1)}}{\|\br^{(k+1)\|_2}}$
\STATE $\bV_{\ell+k+1} = [\bV_{\ell+k}, \bv]$ \Comment{Enlarge the solution subspace} 
\ENDFOR
\STATE \textbf{end function}
\end{algorithmic}
\end{algorithm}

\subsubsection{Solving the Optical Flow Estimation Sub-problem}
To solve the problems \eqref{eq:subproblem_motion_estimation(a)} and \eqref{eq:subproblem_motion_estimation(b)}, we fix $p$ and $q$ and solve the regularized problem with the MMGKS algorithm described in \ref{sec:mmgks_alg}. Since they are similar, we focus on the solution for \eqref{eq:subproblem_motion_estimation(a)}. 

In practice, it is inefficient to solve the combined system \eqref{eq:subproblem_motion_estimation(a)} system for all timesteps at once, instead we solve 
\begin{equation}
    \bs(j)^{(k)} = \argmin_{\mathbf{s} \in R^{2n_s}} \| \bUpsilon^{(k)}(j)\mathbf{s}(j)^{(k)} + \bu^{(k)}_t(j)\|_p^p + \gamma \left\| \begin{bmatrix}\mathbf{L} & \mathbf{L}\end{bmatrix} \bs(j)^{(k)} \right\|_q^q, 
\label{eq:subproblem_motion_estimation_redefined}
\end{equation}
for $t = 1, \ldots n_t-1 $.

We found that for solving a standalone optical flow problem, setting $p = q = 1$ gives optimal results in consistency with results reported in other references \cite{Burger2015OnOF}, however, for the limited angle and single-shot tomography problems we discuss here, choosing $p=q=2$ resulted in more accurate reconstructed image sequence. Hence, for all the experiments we use the later.

The optical flow estimation procedure is summarized in Algorithm \ref{alg: Solve-OF}.

\begin{algorithm}[ht]
\caption{Solve-OF}
\label{alg: Solve-OF}
\begin{algorithmic}[1]
\Require: $\bu(t)$, $\bu(t+1)$, $\bL$
\Ensure: An approximate velocity vector $\bs(t)$ that solves the optical flow equation \eqref{eq:of_lsq}.
\STATE \textbf{function} \textsc{Solve-OF} ({$\bu(t)$, $\bu(t+1)$, $\bL$})
\STATE Compute $\bUpsilon(t)$ and $\bu_t(t)$ 
\STATE $\bs(t)$ = MMGKS ($\bUpsilon(\bu(t))$, $\begin{bmatrix}
 \bL  & \bL
\end{bmatrix}$, $-\bu_t(t)$)
\STATE \Return $\bs(t)$
\STATE \textbf{end function}
\end{algorithmic}
\end{algorithm}

\subsubsection{Regularization Parameter Selection}
\label{sec:regparam}
Both the image reconstruction and the optical flow subproblems require the determination of optimal regularization parameters. The MMGKS algorithm determines this in one of two ways:

\begin{enumerate}
    \item Discrepancy Principle (DP): If an estimate $\delta$ for the norm of the noise in the measurement vector is known, the discrepancy principle computes the regularization parameter for the regularized
 solution of the generic regularization problem \eqref{eq: dynamicEq} by imposing
 \begin{equation}
\mathcal{D}(\bu_{\text{reg}}) := \|\bH \bu_{\text{reg}} - \bb\|_2 = \eta \delta, 
 \end{equation}
 
 where $\eta > 1 (\eta \simeq 1)$ is a safety factor.

\item Generalized Cross-Validation (GCV): The GCV method on the other hand, comes in handy when there is no prior knowledge about the measurement error. The GCV method chooses the regularization parmeter $\lambda$  by minimizing the functional
\end{enumerate}
\begin{equation}
\mathcal{G}(\lambda) :=
 \frac{\|\bH \bu_{\text{reg}}(\lambda) - b\|^2}{\text{trace}(\bI - \bH \bH^\dagger_{\text{reg}}(\lambda) ) 2 }
    \label{eq: gcv_functional}
\end{equation}

In our experiments, for the optical flow estimation subproblem, we select the regularization parameter using GCV, while for the image reconstruction subproblem, we use the discrepancy principle, when information about the noise is known, and  GCV when there is no noise information.

We emphasize here that this automatic and cost-effective selection of the regularization parameter is a key selling point of the MMGKS algorithm in general. 

\subsubsection{Computational Cost}
% \todo{Ask Misha how to structure this.}
The primary computational bottleneck in the MMGKS-OF algorithm is the optical flow computation between consecutive frames within the current solution at each iteration. This cost can be broken down as follows:

\begin{enumerate}
    \item {Optical Flow Equation Solution:} For each of the $n_t - 1$ consecutive image pairs in the current solution's sequence:
    \begin{itemize}
        \item $(n_t - 1) \times 2$ Gradient computations (both spatial and temporal).
        \item $(n_t - 1) \times k_{\text{flow}}$ MMGKS iterations, each consisting of steps 2-4 (detailed below).
    \end{itemize}

    \item {Matrix-Vector Products (at the $k$-th iteration of MM-GKS):}
    \begin{itemize}
        \item One product with matrix $\boldsymbol{\Theta}$.
        \item One product with vector $\bV_{\ell+k}$.
    \end{itemize}

    \item {(Re)orthogonalizations (at the $k$-th iteration of MM-GKS):}
    \begin{itemize}
        \item One orthogonalization for $\bH\bV_{\ell+k}$.
        \item One orthogonalization for $\bP_{\epsilon} \Theta \bV_{\ell+k}$.
        \item $\ell$ orthogonalizations for the initial computation, and one for the subsequent updates, of the $Q_\bH R_\bH$ decomposition.
        \item $\ell$ orthogonalizations for the initial computation, and $\ell+k$ for the subsequent updates, of the $Q_{\mathbf{\Theta}} R_\mathbf{\Theta}$ decomposition.
    \end{itemize}
\item {Least-Squares Solution (per MMGKS iteration):}
    \begin{itemize}
        \item Solution of one reduced least-squares problem.
    \end{itemize}
\end{enumerate}

\begin{table}[h]
    \centering
    \begin{tabular}{|c|c|}
        \hline
        Notation & Meaning \\ \hline
        $\ell$ & Dimension of initial subspace \\ \hline
        $k$ & Iteration number \\ \hline
        $k_{\text{of}}$ & Number of iterations of the optical flow solver \\ \hline
        $n_x, n_y$ & Image dimensions\\ \hline
        $n_s$ & $n_xn_y$ \\ \hline
        $n_t$ & Number of images in the sequence\\
        \hline
        $n$ & $n_sn_t$ \\ \hline
        $\bu_t(t)$ & Derivative of $\bu$ w.r.t. time \\ \hline
         $\bu_t$ & $\begin{pmatrix} \bu_t(1) & \cdots & \bu_t(n_t-1)\end{pmatrix}^\top$\\ \hline
        $\bs(t)$ & The velocity vector for all the pixels in $\bu(t)$ at time $t$\\ \hline
        $\bs$  & $\begin{pmatrix} \bs(1) & \cdots & \bs(n_t-1)\end{pmatrix}^\mathrm{T}$ \\ \hline
        $\bu(t)$ & A single vectorized 2D image \\ \hline
        $\bu$ & $\text{vec}([\bu{(1)}, \ldots, \bu{(n_t)}])$ \\ \hline   \end{tabular}
    \caption{Table of Notations}
    \label{tab:notation_table}
\end{table}
\subsubsection{Strategies for Reducing Computational Costs}
\label{sec: comp_cost_reduction}

Despite the efficiency of the MMGKS-OF solver, which reduces the solution space to a lower dimensional one, solving the optical flow equations for the velocity vectors at each time step can still be computationally expensive. To address this challenge, we implement the following strategies to reduce computational cost while maintaining high algorithmic performance:
\begin{itemize}

\item {Reduced Frequency of Optical Flow Computation}:
It is not necessary to compute the optical flow at every iteration of the MMGKS-OF algorithm. Instead, we introduce a parameter $\tau$, which controls the frequency of optical flow updates. This reduces the overall computational effort by skipping redundant calculations during iterations where the solution does not significantly change.

\item {Reverse Optical Flow Approximation}:  
To further reduce computational costs, we compute the optical flow in only one direction, say \(\mathbf{s}\), and then use an approximation, as described in \eqref{eq:s_prime_from_s} to estimate \(\mathbf{s}'\). This strategy effectively halves the number of flow computations without compromising the solution's accuracy.

\item {Parallelized Optical Flow Computation}:  
For iterations where the optical flow is computed, we parallelize the process across image pairs in the sequence. This is achieved by dividing the workload among multiple processors or GPUs, significantly accelerating the computations, especially for long image sequences or high-resolution data.

\item {Rescaling Large Images for Optical Flow Computation}:  
For large images, depending on a scaling parameter $\alpha$, SOLVE-OF first rescales them to a smaller resolution before computing the optical flow. This step reduces the computational burden significantly, as the size of the system to be solved scales with the image dimensions. After obtaining the optical flow on the smaller images, the results is upscaled to match the original resolution.

\end{itemize}

\section{Numerical Experiments}
\label{sec:experiments}
To illustrate the performance of the MM-GKS-OF algorithm, we consider four different variations shown in Table \ref{tab: methods} based on the choice of the regularization operator $\Psi$ and compare their performance with non-OF methods on applications in computerized tomography including
limited angle tomography and single shot tomography. 

The notation for these variations are discussed further in \cref{rem:of-variations}

\begin{remark}
\label{rem:of-variations}
    In \cref{tab: methods}, $\T{U}$ represents the tensor representation of $\bu$ as described in \cite{pasha2023computational}, $\bL_{h}$, $\bL_{v}$, $\bL_{t}$ are discretizations of the first derivatives in the vertical, and time directions respectively as defined in \cref{eq: L}. \emph{Mode-$j$ product} that defines the operation of multiplying a tensor $\T{U}$ (of appropriate dimensions) by a matrix $\bL_j \in \R^{r\times n_j}$ for $j=1,2,3$ is denoted by $\times_j$.
\end{remark}

\begin{table}
    \centering
    \begin{tabular}{|c|c|c|}
    \hline
        Method & Notation & Regularization term \\ \hline
        MM-GKS & M & $ \boldsymbol{\Psi} \bu$ \\ \hline
        MM-GKS-OF & M-OF & $\boldsymbol{\Psi} \bu +  \mathcal{M}(\bu) $\\ \hline
        AnisoTV & D$_1$ & $ \mathcal{R}_1(\bu) := \|\T{U} \times_1 \bL_{v}\|_1 + \|\T{U} \times_2 \bL_{h}\|_1 + \|\T{U} \times_3 \bL_{t}\|_1$ \\ \hline
        AnisoTV-OF& D$_1$-OF & $\mathcal{R}_1(\bu)  + \mathcal{M}(\bu)$ \\ \hline
        IsoTV & D$_2$ &  $\mathcal{R}_2(\bu) := \sum_{i=1}^{n_s} \sum_{j=1}^{n_t} \|\T{Y}_{i,j,:}\|_2 +\|\T{U}\times_3\bL_t\|_1$\\ \hline
        IsoTV-OF& D$_2$-OF & $\mathcal{R}_2(\bu)  + \mathcal{M}(\bu)$\\ \hline
        GS & D$_3$ & $\mathcal{R}_3(\bu) := \|\bL_s\bU\|_{2,1}$ \\ \hline
        GS-OF & D$_3$-OF & $\mathcal{R}_3(\bu)  + \mathcal{M}(\bu)$\\ \hline
    \end{tabular}
    \caption{Summary of methods used in the numerical experiments.}
    \label{tab: methods}
\end{table}

When a ground truth is available, the regularization parameter is determined using the discrepancy principle; otherwise, generalized cross-validation (GCV) is employed. All methods were executed for a maximum of 200 iterations. To manage computational resources (as detailed in Section \ref{sec: comp_cost_reduction}), optical flow (OF) calculations were performed only every 20 iterations. Furthermore, prior to OF computation, all images were downsampled by a factor of 4 or 2 to achieve dimensions smaller than 50x50 pixels. The resulting velocity vectors were subsequently upscaled to their original dimensions. While these optimizations slightly impact OF performance, the effect is minimal, as these methods still outperform the alternatives.

In cases where the forward operator and data were synthetically generated, we mitigated the inverse crime by adding simulated measurement noise with a magnitude of $10^{-1}$ to the simulated sinograms.

Except for Test 5, all experiments were conducted on a Dell XPS 16 with the following specifications: Processor: Intel Core i9-12900HK (14 cores, 20 threads, 2.5 GHz base, up to 5.0 GHz boost); Memory: 64 GB DDR5 RAM; Graphics Card: NVIDIA GeForce RTX 3060, 6 GB GDDR6.

\paragraph{Discussion on the selection of the numerical examples}
Two classes of examples are considered in this study: limited-angle tomography and single-shot tomography. For limited-angle tomography, three tests are conducted: two using synthetic data and one using real data. The first synthetic test (Test 1) investigates the impact of varying the number of projection angles on the performance difference between MMGKS-OF and other non-OF methods, while also evaluating the quality of the optical flow estimated by MMGKS-OF. The second synthetic test (Test 2) uses a real image as a phantom to examine how the selection of projection angles for the tomographic forward operator influences MMGKS-OF performance. The third test (Test 3) employs a real dataset with sinograms obtained from actual measurements. Although the ground truth image is unknown, allowing only for qualitative assessment, this example provides valuable insight into the proposed algorithm's performance in real-world scenarios with noisy data. For single-shot tomography, two synthetic examples are considered. The first (Test 4) uses a dataset obtained from publicly available code accompanying a paper on dynamic shape sensing. While the underlying method differs, this dataset provides a valuable test case for MMGKS-OF performance on a different type of data. The second synthetic example (Test 5) uses data derived from a paper on joint image reconstruction, which employs a similar framework but a different solution algorithm. The primary aim of using this dataset is for comparison, particularly regarding sensitivity to the regularization parameter.

\paragraph{Quality measures and stopping criteria} 
To evaluate the quality of the reconstructed images, we use two metrics: the relative reconstruction error (RRE) and the structural similarity index (SSIM). The RRE between the recovered $k$-th iterate $\bu^{(k)}$ and the ground truth image $\bu_{\text{true}}$ is defined as 
\begin{equation}
    \text{RRE}(\bu^{(k)}, \bu_{\text{true}}) = \frac{\|\bu^{(k)} - \bu_{\text{true}}\|}{\|\bu_{\text{true}}\|_2}
\end{equation}
SSIM evaluates similarity by comparing local patterns of pixel intensities, considering luminance, contrast, and structure. Unlike metrics based on pixel-wise differences, SSIM focuses on changes in structural information. SSIM values range from 0 to 1, with 1 indicating perfect similarity and 0 indicating no similarity. The detailed SSIM formula is provided in \cite{SSIM}.

Iterations are halted after a maximum of 200 iterations. To ensure fair comparison, the initial subspace dimension in the MMGKS-OF algorithm is set to 10 for all experiments, and the MMGKS algorithm in \textsc {SOLVE-OF} is executed for 20 iterations per call.

All test problems and algorithm implementations in Python will be made publicly available at \url{https://github.com/mpasha3/Dynamic_Inverse_Problems_with_Optical_Flow.git} once the manuscript is accepted to the journal.

\subsection{Computerized Limited Angle Tomography}
\subsubsection{Test 1: Simulated Moving Blocks Dataset}

We generate a 2D dynamic phantom with four blocks moving at constant speed within a bounded space for 12 frames, two of the blocks move faster than the other two. We initialize the positions of the blocks within a frame and simulate rigid motion of the blocks in different directions and at different speeds within the box. The resulting phantom is a 3D array with dimensions of $90 \times 90 \times 12$.

\begin{figure}[htbp]
    \centering
    \includegraphics[width=0.9\linewidth]{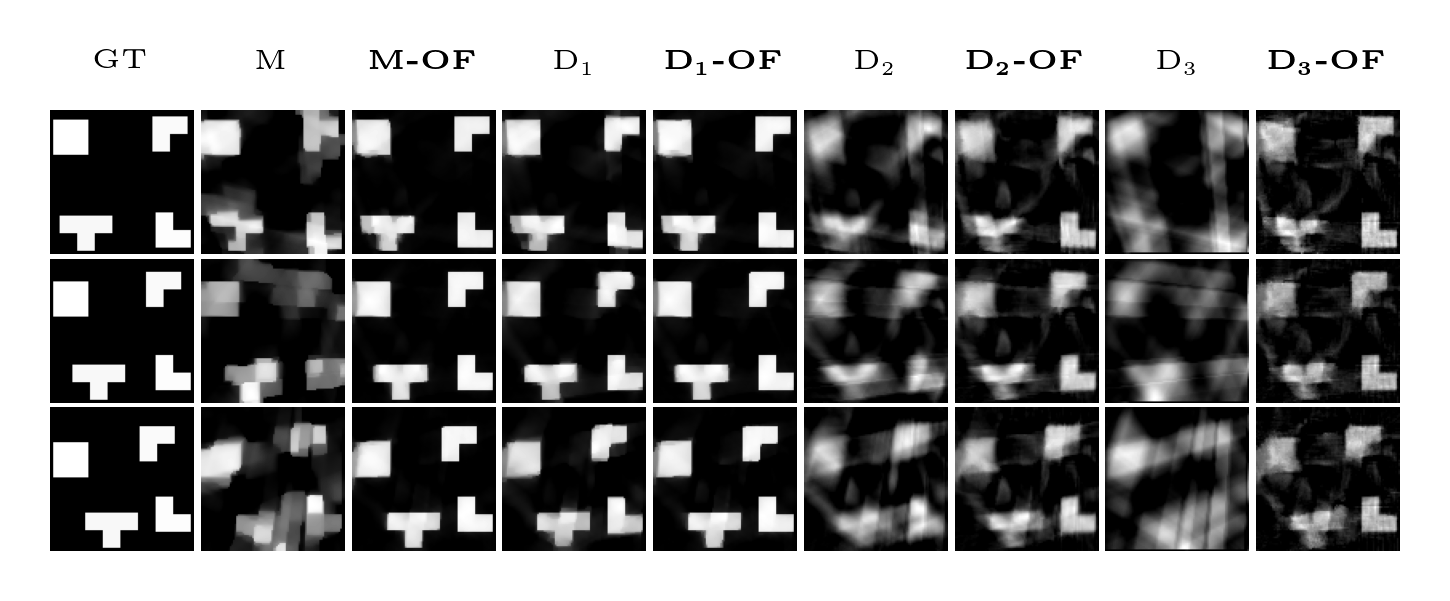}
    \caption{Test 1: Reconstructed images  at $t = \{2,6,10\}$ when $n_{\text{views}} = 3$.}
    \label{fig:test1-rec_images}
\end{figure}

\paragraph{Forward Model and Projection Geometry}
To generate the forward operator at each timestep, we use the Trips-Py toolbox \cite{pashapythonpackage} to create a 2D flat fan beam geometry with 117 rays departing from each angle in a set of projection angles. For the projection angles, we set $\zeta := \frac{180}{n_{\text{views}}}$ and define the interval $\Theta_t := \left[\frac{\zeta(t - 1)}{n_t}, \frac{\zeta(t - 1)}{n_t} + 180^\circ\right)$. Each interval has a width of $180^\circ$ and is shifted by $\frac{\zeta}{n_t}$ from the previous interval. For example, when the number of projection angles, $n_{\text{views}} = 3$ and $n_{t} = 12$, $\zeta = \frac{180}{3} = 60$, and $\Theta_t := \left[\frac{60(t - 1)}{12}, \frac{60(t - 1)}{12} + 180^\circ\right)$, so $\Theta_1 = [0,180^\circ)$, $\Theta_2 = [5,185^\circ)$, ..., $\Theta_{12} = [55,235^\circ)$. For given $n_{\text{views}}$, we obtain the forward problem $\bH^{t}\bu{(t)} + \gamma^{t} = \bb^{t}$, $t = 1,2,\ldots,12$, where $\bH^{t} \in \R^{117 \cdot n_{\text{views}} \times 8100}$, $\bb^{t} \in \R^{2,170} $. The block forward operator $\bH$ has size $1,524 \cdot n_{\text{views}} \times 97,200$, and the measurement vector $\bb$ has size 97,200. 
\paragraph{Mitigating Inverse Crime} In order to mitigate the risk of inverse crime, we employed several strategies in our reconstruction process. Inverse crime occurs when the same model is used for both generating synthetic data and reconstructing it, leading to overly optimistic results.
First, we used different projector models to create the projection geometry and for the actual reconstruction. This ensured that the forward model and the reconstruction model were not identical.
In addition, we introduced a slight misalignment in the projection geometry. By using different angles for the misaligned and standard projection geometries, we simulated more realistic conditions and avoided the perfect alignment that can lead to inverse crime.
\paragraph{Results and Discussion}
We vary the number of projection angles $n_{\text{views}}$ and inspect the performance of the OF algorithms against their non-OF counterparts. Figure \ref{fig:test1-rec_images} shows the reconstructed images obtained in selected time steps when $n_{\text{views}} = 3$. Table \ref{tab:test1-rre_and_ssim} presents the mean RREs and SSIM indices achieved by the methods. Our results consistently demonstrate the superior performance of optical flow (OF) methods compared to their non-OF counterparts, with MMGKS-OF achieving the best overall results. We also observe a diminishing return in the improvement offered by MMGKS-OF over non-OF methods as the amount of available information increases (i.e., with a higher number of projection angles). This is because with more data, the additional information provided by the motion regularization term becomes somewhat redundant, as all methods converge more quickly. Nevertheless, the MMGKS method proves particularly advantageous in severely limited-angle scenarios, as evidenced by the substantial performance difference observed when $n_{\text{views}} = 5$.
\begin{figure}[t]
    \centering
    \begin{subfigure}{0.45\textwidth}
        \centering
        \includegraphics[width=\linewidth]{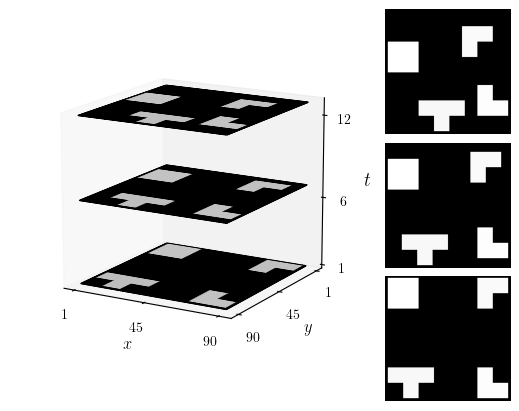}
        \caption{Test 1}
    \end{subfigure}%
    \hspace{0.01\textwidth} % Adjust the space between figures
    \begin{subfigure}{0.45\textwidth}
        \centering
        \includegraphics[width=\linewidth]{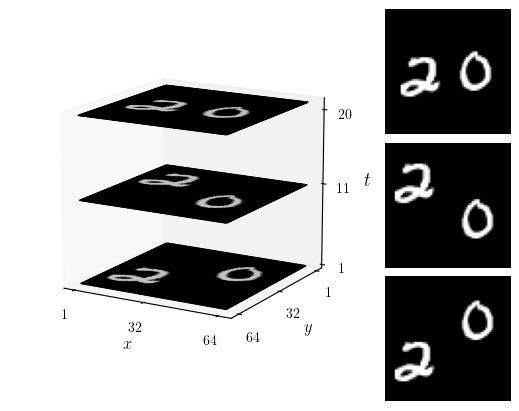}
        \caption{Test 2}
    \end{subfigure}%
    \caption{Slices of ground truth images of (a) Test 1 at $t = \{1,6,12\}$, (b)Test 2at $t = \{1, 11, 20\}$}
    \label{fig:ground_truth_slices_tests_1_and_2}
\end{figure}

\begin{figure}[htbp]
    \centering
    \begin{subfigure}{0.3\textwidth}
        \centering
        \includegraphics[width=\linewidth]{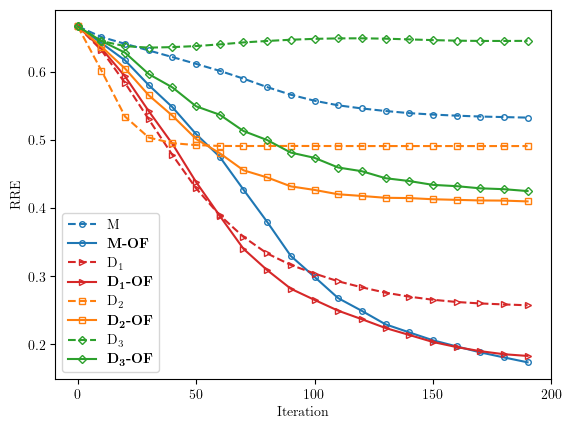}
        \caption{$n_\text{views} = 3$}
    \end{subfigure}%
    \hspace{0.01\textwidth} 
    \begin{subfigure}{0.3\textwidth}
        \centering
        \includegraphics[width=\linewidth]{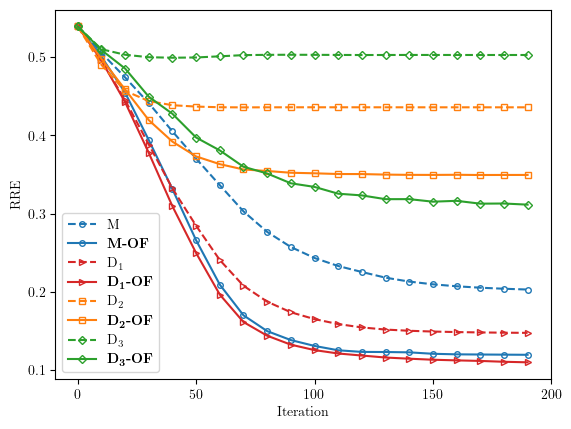}
        \caption{$n_\text{views} = 5$}
    \end{subfigure}%
    \hspace{0.01\textwidth} 
    \begin{subfigure}{0.3\textwidth}
        \centering
        \includegraphics[width=\linewidth]{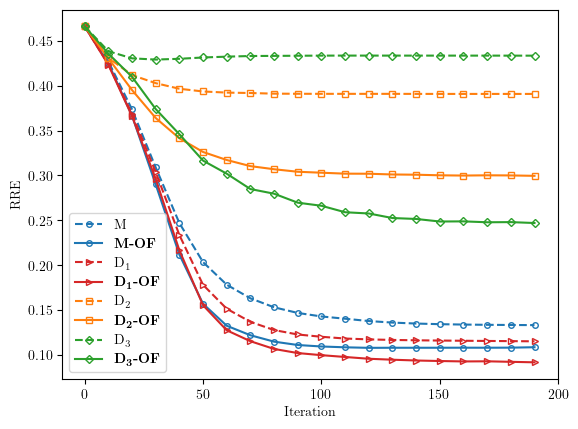}
        \caption{$n_\text{views} = 7$}
    \end{subfigure}
    \caption{Test 1: RRE convergence plots}
    \label{fig:test1-rre_plots}
    \end{figure}

\begin{figure}[htbp]
    \centering
    \begin{subfigure}{0.45\textwidth}
        \centering
        \includegraphics[width = \textwidth]{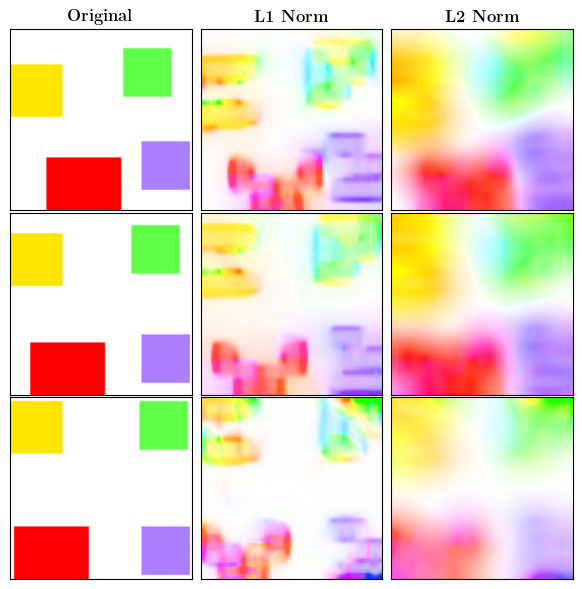}
        \caption{Optical flow}
    \end{subfigure}%
    \hspace{0.01\textwidth}  
    \begin{subfigure}{0.45\textwidth}
        \centering
        \includegraphics[width = \textwidth]{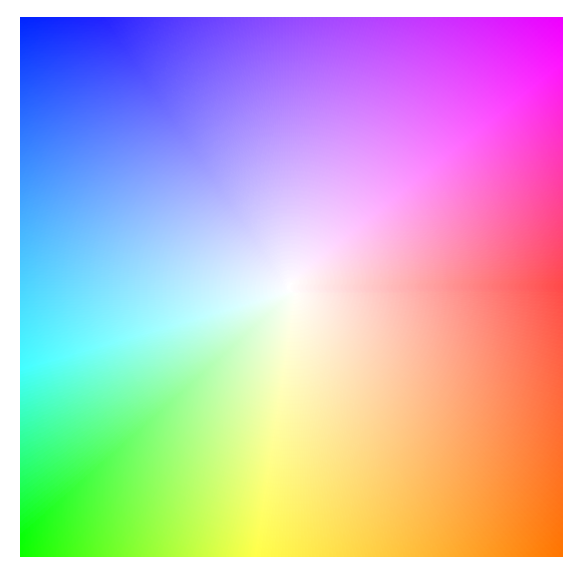}
        \caption{Optical flow color bar}
    \end{subfigure}
    \caption{Test 1: Estimated optical flow when $n_\text{views} = 3$ at $t = 1, 5, 10$.}
    \label{fig:test1-of_plots}
\end{figure}

\begin{table}[t]
    \centering
    \begin{tabular}{|c|c|c|c|c|c|c|c|c|}
\hline
\multirow{ 2}{*}{$n_{\text{views}}$} & \multicolumn{8}{|c|}{RRE} \\
\cline{2-9}
& M & \textbf{M-OF}& D$_1$&\textbf{D$_1$-O}F& D$_2$& \textbf{D$_2$-OF}& D$_3$& \textbf{D$_3$-OF}\\
\hline
3 & 0.528 & 0.167  & 0.251  & 0.177 & 0.491 & 0.409 & 0.644 & 0.423 \\
\hline
5 & 0.590 & 0.341  & 0.376  & 0.319 & 0.430 & 0.376 & 0.607 & 0.405 \\
\hline
7 & 0.582 & 0.330 & 0.370 & 0.313 & 0.424  & 0.371 & 0.595 & 0.402 \\
\hline
 & \multicolumn{8}{|c|}{SSIM} \\

\hline
3 & 0.608 & 0.358  & 0.388  & 0.336 & 0.443 & 0.398 & 0.629 & 0.433 \\
\hline
5 & 0.590 & 0.341  & 0.376  & 0.319 & 0.430 & 0.376 & 0.607 & 0.405 \\
\hline
7 & 0.582 & 0.330 & 0.370 & 0.313 & 0.424  & 0.371 & 0.595 & 0.402 \\
\hline
\end{tabular}
\caption{\textit{\textbf{Test 1}: Mean RRE and SSIM for 3, 5 and 7 projection angles.}}
    \label{tab:test1-rre_and_ssim}
\end{table}

\subsubsection{Test 2: MNIST Moving Numbers Dataset} In this example we consider the MNIST dataset found in \cite{moving_mnist} which consists of 10,000 video sequences of length 20. Each frame in the sequence is of resolution 64×64 pixels and contains two moving digits. From this dataset, we select a single sequence as our ground truth.

\paragraph{Forward Model and Projection Geometry}
To simulate the forward projection process, we employ the Trips-Py toolbox \cite{pashapythonpackage} to construct a 2D flat fan-beam geometry. We consider three different scenarios for the selection of projection angles:  

\begin{enumerate}

\item Fixed Angle Set: All intervals are the same for each timestep, i.e. $\Theta_t = [0,180)$ for t = $1,\ldots n_t$. 

    \item Equal Intervals: The interval $[0,180)$ is split into $n_t$ equal intervals each timestep is assigned one of these sub-intervals, i.e. $\Theta_t = [\frac{180 (t-1)}{n_t},\frac{180 t}{n_t})$ for t = $1,\ldots n_t$, e.g. for $n_{\text{views}} = 4$, and $n_{t} = 3$, $\Theta_1 = [0,60)$,$\Theta_2 = [60,120)$, and $\Theta_3 = [120 \degree,180)$.

    \item Shifted Intervals: We set $\zeta: = \frac{180}{n_{views}}$ and, we define the interval $\Theta_t := \left[\tfrac{\zeta(t - 1)}{n_t}, \tfrac{\zeta(t - 1)}{n_t} + 180 \degree \right)$, i.e., each interval has a width of $180 \degree$ and is a shift by $\tfrac{\zeta}{n_t}$ from the previous interval, e.g. for $n_{\text{views}} = 4$, and $n_{t} = 3$, $\zeta = \frac{180}{4} = 45$,  and $\Theta_t := [15(t - 1), 15 (t-1) + 180)$, so $\Theta_1 = [0,180 \degree)$,$\Theta_2 = [15,195 \degree)$, and $\Theta_3 = [30 \degree, 210\degree)$.
\end{enumerate}
    \paragraph{Results and Discussion}
We evaluate the performance of all eight methods under each of the three angle selection scenarios. In all scenarios, we observe that the OF methods consistently outperform their non-OF counterparts. However, the most significant improvement between OF and non-OF methods is seen in scenario (b), where the projection angles are binned. This significant improvement can likely be attributed to the fact that binning reduces the effective angular coverage, introducing sparsity and gaps in the projection data. Non-OF methods struggle in such scenarios due to their reliance solely on the provided projections, which limits their ability to recover missing information. In contrast, OF methods effectively leverage motion estimation to exploit structural relationships and temporal coherence, compensating for the degraded angular resolution and providing better reconstructions. The ability of OF methods to fill in these gaps aligns particularly well with the challenges posed by binned data, making this the scenario where their advantages are most pronounced.

Interestingly, in this test, the M-OF method without the temporal derivative regularizer outperforms all the D$_i$-OF methods. This is likely due to the larger movement between subsequent frames, which can make the spatio-temporal derivative regularizer less accurate in the temporal dimension. This suggests that while the D$_i$-OF methods are generally effective, the M-OF method with the spatial derivative regularizer may be more suitable for scenarios with significant motion. 

\begin{figure}[htbp]
    \centering
    \begin{subfigure}{0.3\textwidth}
        \centering
        \includegraphics[width=\linewidth]{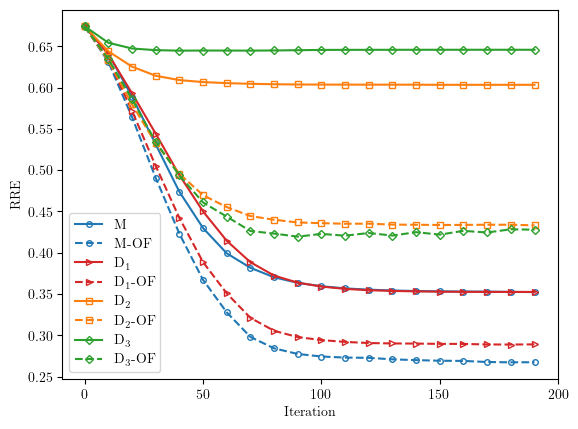}
        \caption{Case a}
    \end{subfigure}%
    \hspace{0.01\textwidth} 
    \begin{subfigure}{0.3\textwidth}
        \centering
        \includegraphics[width=\linewidth]{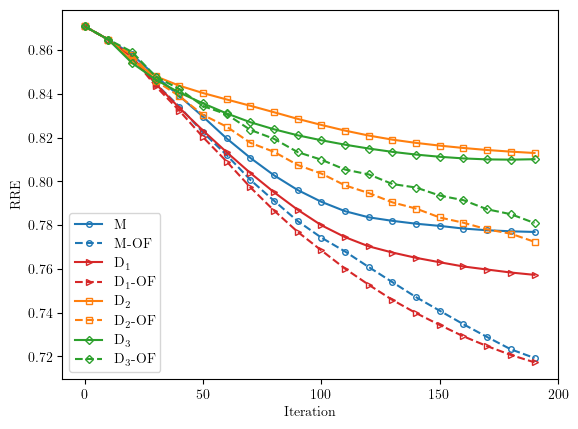}
        \caption{Case b}
    \end{subfigure}%
    \hspace{0.01\textwidth} 
    \begin{subfigure}{0.3\textwidth}
        \centering
        \includegraphics[width=\linewidth]{figures_/example_2_mnist/test2_case_b.png}
        \caption{Case c}
    \end{subfigure}
    \caption{Test 2: RRE convergence plots for MMGKS, MMGKS-OF, ANISO, ANISO-OF, ISO, ISO-OF, GS and GS-OF when the angles are selected according to scenarios a, b and c.}
    \label{fig:test2-rre_plots}
    \end{figure}

\begin{figure}[t]
    \centering
    \includegraphics[width=0.9\linewidth]{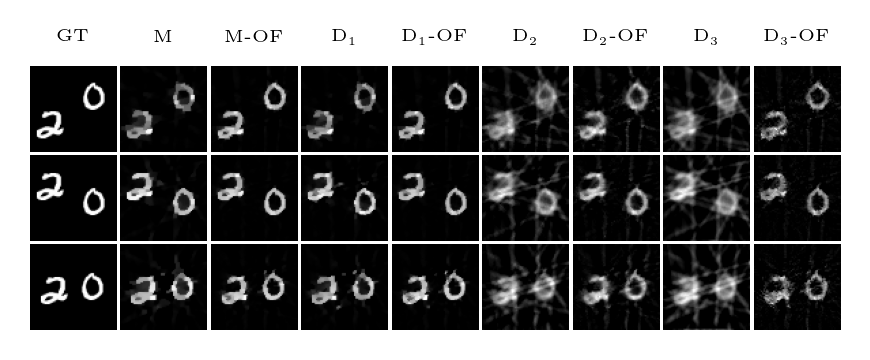}
    \caption{Test 2: Reconstructed images  at $t = \{2,11,18\}$ when angles are selected according to Case a.}
    \label{fig:test2-rec_images}
\end{figure}

Figure \ref{fig:test2-rre_plots} illustrates the convergence of the RRE in all three angle selection scenarios, and \ref{fig:test2-rec_images} shows the reconstructed images at selected time steps, while Table \ref{tab:test2-rre_and_ssim} shows the mean RREs and SSIM indices achieved by the methods.

\begin{table}[t]
    \centering
    \begin{tabular}{|c|c|c|c|c|c|c|c|c|}
\hline
\multirow{2}{*}{Case} & \multicolumn{8}{|c|}{RRE} \\
\cline{2-9}
& M & \textbf{M-OF} & D$_1$ & \textbf{D$_1$-OF} & D$_2$ & \textbf{D$_2$-OF} & D$_3$ & \textbf{D$_3$-OF} \\
\hline
a & 0.351 & 0.267 & 0.351 & 0.288 & 0.603 & 0.433 & 0.646 & 0.437 \\
\hline
b & 0.776 & 0.714 & 0.756 & 0.713 & 0.812 & 0.768 & 0.810 & 0.778 \\
\hline
c & 0.384 & 0.293 & 0.373 & 0.307 & 0.615 & 0.440 & 0.654 & 0.434 \\
\hline
 & \multicolumn{8}{|c|}{SSIM} \\
\hline
a & 0.920 & 0.957 & 0.922 & 0.949 & 0.700 & 0.863 & 0.644 & 0.868 \\
\hline
b & 0.465 & 0.572 & 0.508 & 0.577 & 0.406 & 0.490 & 0.401 & 0.493 \\
\hline
c & 0.902 & 0.948 & 0.910 & 0.942 & 0.688 & 0.859 & 0.635 & 0.869 \\
\hline
\end{tabular}
\caption{\textit{\textbf{Test 2}: Mean RRE and SSIM for Cases a, b, and c.}}
\label{tab:test2-rre_and_ssim}
\end{table}

\subsubsection{Test 3: Emoji Dataset}  
\begin{figure}[t]
    \centering
    \includegraphics[width=0.9\linewidth]{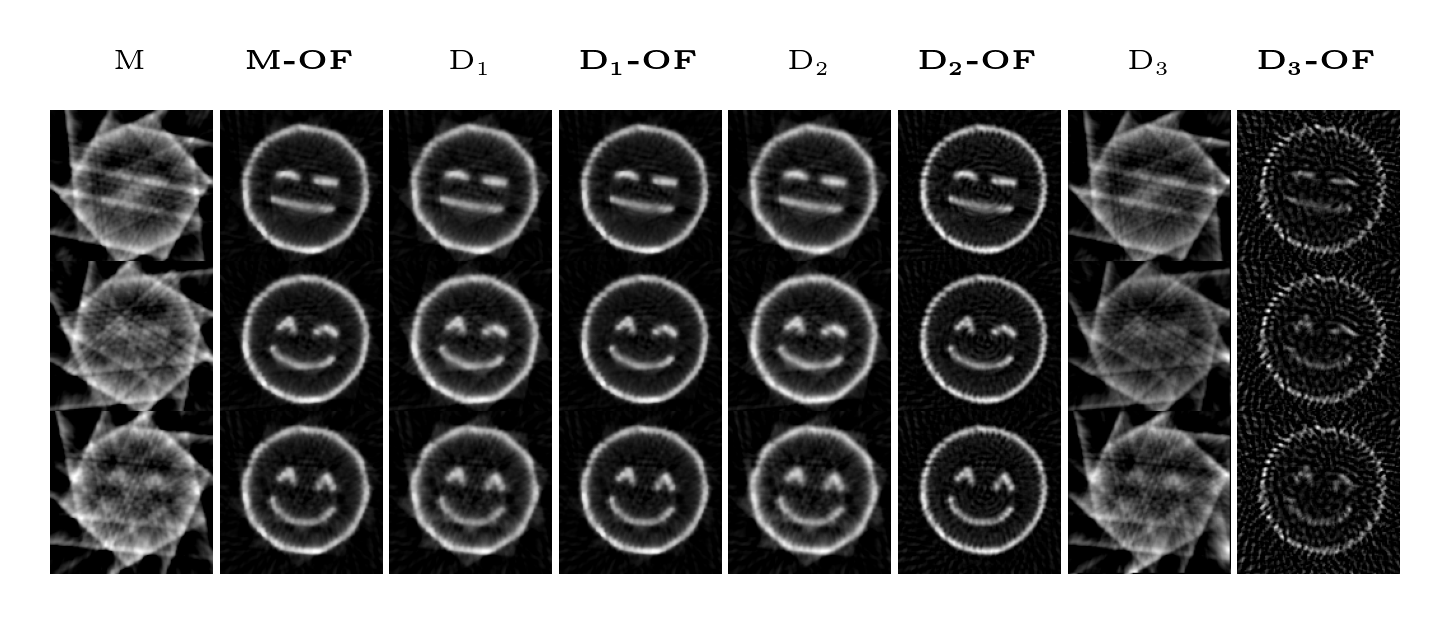}
    \caption{Test 3: Reconstructed Images  at $t = \{2,17,31\}$}
    \label{fig:test3-rec_images}
\end{figure}

This example considers \emph{real data} of an ``emoji" phantom measured at the University of Helsinki [42].
\paragraph{Forward Model and Projection Geometry}
We modify the data to determine a limited angle problem by limiting the number of projection angles to 10. From the dataset ``DataDynamic $128 \times 30$.mat", we generate the
 problems $\bH^{t}\bu{(t)} + \gamma^{t} = \bb^{t}$, $t = 1,2,\ldots,33$ where $\bH^{t} \in \R^{2,170 \times 16,384}$, $\bb^{t} \in \R^{2,170} $ are defined by taking 217
 fan-beam projections around only 10 equidistant angles in $[0,2 \pi)$. The block forward operator $\bH$ has size $71,610 \times 540,672$, and the measurement vector $\bb$ has size 71,610. 
 \paragraph{Results and Discussion}
 There is no ground truth available for the phantom, so we simply present the reconstructed images in \ref{fig:test3-rec_images}. However, from observation it is evident that the OF methods again outperform the others, with the ISO-OF method giving the best results this time.

\subsection{Computerized Single-shot Tomography}

\subsubsection{Test 4:  Single-shot tomography of synthetic rigid motion }
\begin{figure}[t]
    \centering
    \includegraphics[width=0.9\linewidth]{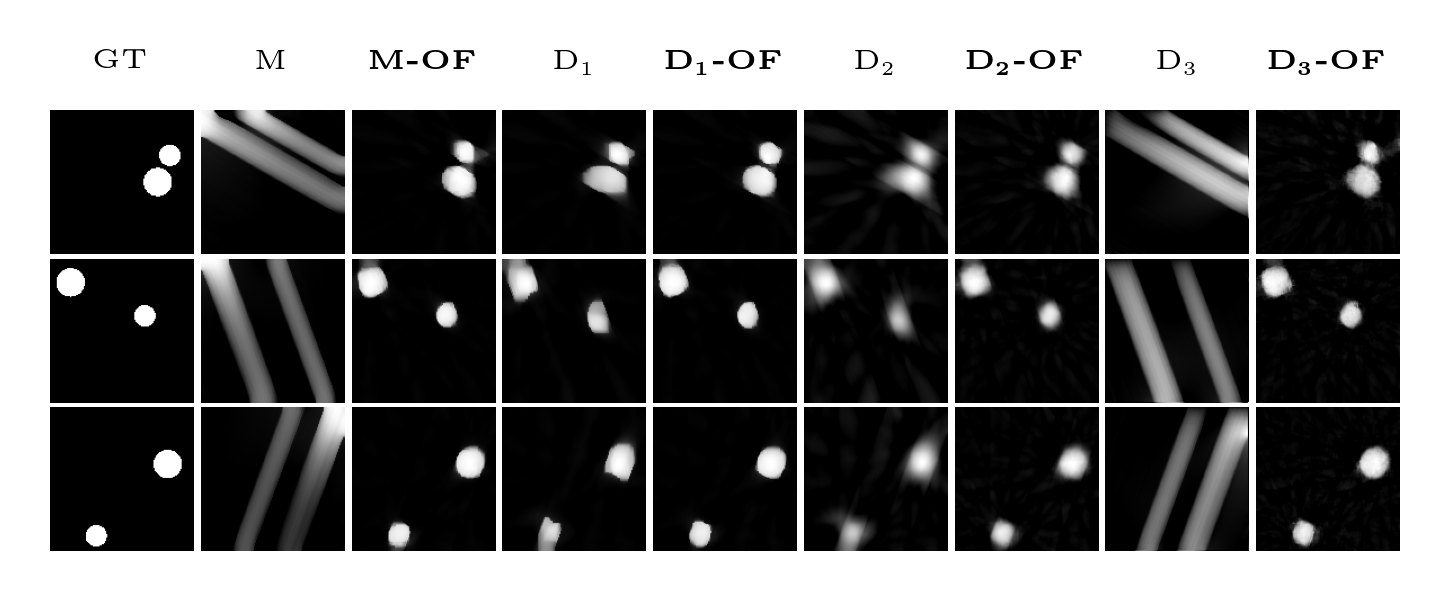}
    \caption{Test 4: Reconstructed Images at  at $t = \{30,64,98\}$}
    \label{fig:test4-rec_images}
\end{figure}

\begin{figure}[htbp]
    \centering
    \begin{subfigure}{0.45\textwidth}
        \centering
        \includegraphics[width=\linewidth]{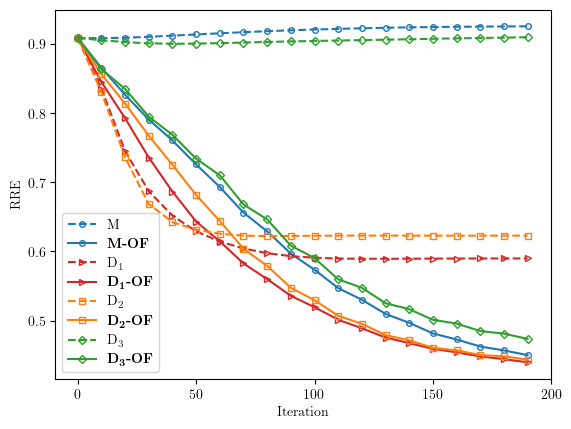}
        \caption{Test 4}
        \label{fig:test4-rre_plots}
    \end{subfigure}%
    \hspace{0.01\textwidth} 
    \begin{subfigure}{0.45\textwidth}
        \centering
        \includegraphics[width=\linewidth]{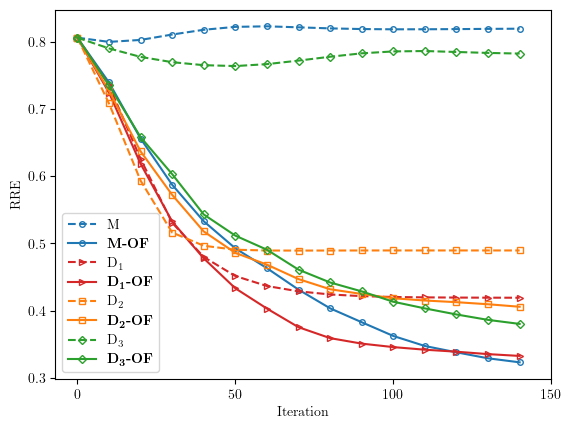}
        \caption{Test 5}
        \label{fig:test5-rre_plots}
    \end{subfigure}%
    \caption{Test 4 (L) and Test 5 (R): RRE convergence plots.}
    \end{figure}

For this example, we consider the synthetic rigid motion dataset generated and used in \cite{figueiredo2003algorithm}.  
\paragraph{Forward Model and Projection Geometry}
Using the code supplied in the original paper, we generate a sequence of 128 frames of 128 by 128 images, the forward operators and the sinograms.
\paragraph{Results and Discussion}
\cref{fig:test4-rec_images} shows the reconstructed images obtained in selected time steps. Table \ref{tab:test4-rre_and_ssim} presents the mean RREs and SSIM indices achieved by the methods and Figure \ref{fig:test4-rre_plots} shows the convergence of the RRE. Our results demonstrate that again, the OF methods consistently outperform their non-OF counterparts, with the ANISO-OF method achieving the best performance.

\begin{figure}[t]
    \centering
    \begin{subfigure}{0.45\textwidth}
        \centering
        \includegraphics[width=\linewidth]{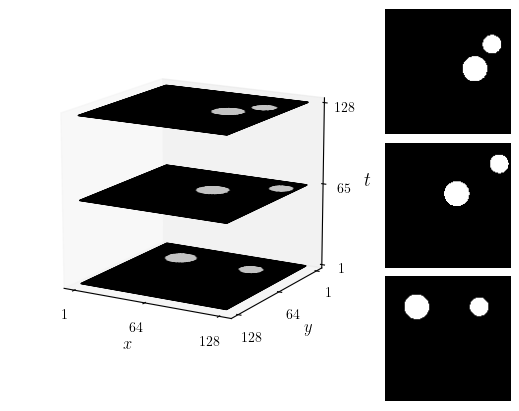}
        \caption{Test 4}
    \end{subfigure}%
    \hspace{0.01\textwidth} % Adjust the space between figures
    \begin{subfigure}{0.45\textwidth}
        \centering
        \includegraphics[width=\linewidth]{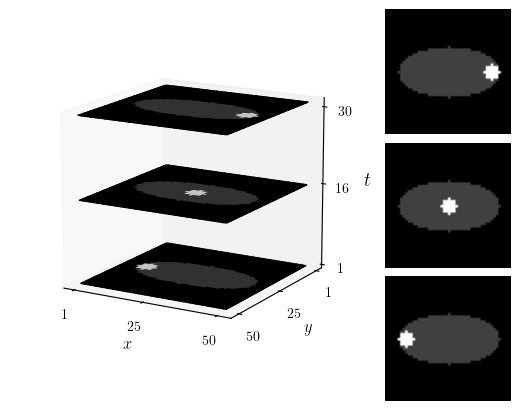}
        \caption{Test 5}
    \end{subfigure}%
    \caption{Slices of ground truth images of (a) Test 4 at $t = \{1,65,128\}$, (b) Test 5 at $t = \{1, 16, 30\}$}
    \label{fig:ground_truth_slices_tests_4_and_5}
\end{figure}

\begin{table}[t]
    \centering
    \begin{tabular}{|c|c|c|c|c|c|c|c|c|}
\hline
& M & \textbf{M-OF} & D$_1$ & \textbf{D$_1$-OF} & D$_2$ & \textbf{D$_2$-OF} & D$_3$ & \textbf{D$_3$-OF} \\
\hline
RRE & 0.925 & 0.439 & 0.589 & 0.433 & 0.622 & 0.439 & 0.909 & 0.461 \\
\hline
SSIM & 0.212 & 0.877 & 0.758 & 0.881 & 0.698 & 0.864 & 0.239 & 0.848 \\
\hline
\end{tabular}
\caption{\textit{\textbf{Test 4}: Mean RRE and SSIM.}}
\label{tab:test4-rre_and_ssim}
\end{table}

\subsubsection{Test 5: Single-shot tomography of pinball dataset}
We simulate a pinball CT dataset similar to the one used in \cite{burger2017variational} consists of forward operators and sinograms of a phantom, a 50 by 50 image of a uniform and rigid ball in a stationary ellipse moving from the left side of the image frame to the right side. Figure \ref{fig:ground_truth_slices_tests_4_and_5}(a) shows the ground-truth images at a selection of time steps. 

\paragraph{Forward Model and Projection Geometry}
From the dataset provided, we generate the problems  $\bH^{t}\bu{(t)} + \gamma^{t} = \bb^{t}$, $t = 1,2,\ldots,33$ where $\bH^{t} \in \R^{75 \times 2500}$, $\bb^{t} \in \R^{75} $ are defined by taking projections around only a random single angle in $[0,2\pi)$. The block forward operator $\bH$ has size $2,250 \times 75,000$, and $\bb$ has size $2,250$.

\begin{figure}[t]
    \centering
    \includegraphics[width=0.9\linewidth]{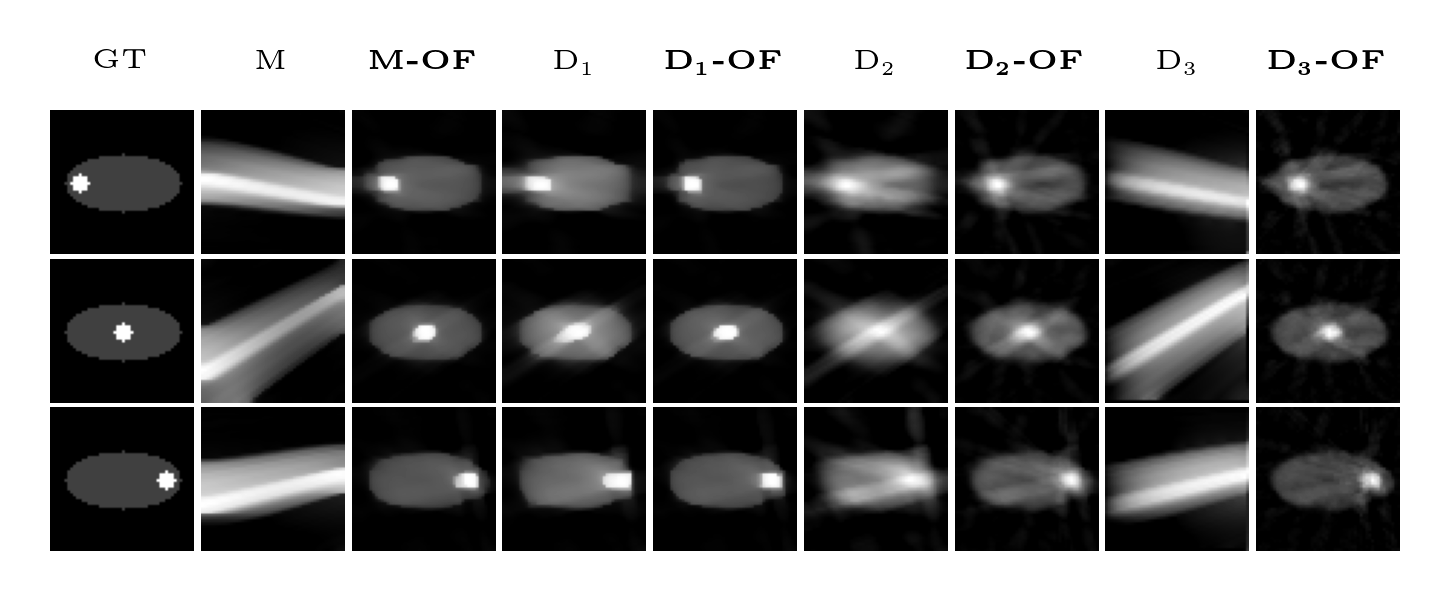}
    \caption{Test 5: Reconstructed Images at $t = 2, 14, 28$.}
    \label{fig:test5-rec_images}
\end{figure}

\paragraph{Results and Discussion}
Figure \ref{fig:test5-rec_images} shows the reconstructed images obtained in selected time steps. Table \ref{tab:test5-rre_and_ssim} presents the mean RREs and SSIM indices achieved by the methods and Figure \ref{fig:test5-rre_plots} shows the convergence of the RRE. Our results demonstrate that again, the OF methods consistently outperform their non-OF counterparts, with the M-OF method achieving the best performance. In the appendix, we also show the performance of the method proposed in \cite{burger2017variational} , highlighting the sensitive dependence of their work on a parameter in the problem formulation.

\begin{table}[t]
    \centering
    \begin{tabular}{|c|c|c|c|c|c|c|c|c|}
    \hline
    & M & \textbf{M-OF} & D$_1$ & \textbf{D$_1$-OF} & D$_2$ & \textbf{D$_2$-OF} & D$_3$ & \textbf{D$_3$-OF} \\
    \hline
    \text{RRE} & 0.815 & 0.314 & 0.418 & 0.323 & 0.489 & 0.399 & 0.777 & 0.368 \\
    \hline
    \text{SSIM} & 0.301 & 0.926 & 0.863 & 0.922 & 0.806 & 0.880 & 0.384 & 0.896  \\
    \hline
    \end{tabular}
    \caption{\textit{\textbf{Test 5}: Mean RRE and SSIM.}}
    \label{tab:test5-rre_and_ssim}
\end{table}

\section{Conclusions and Outlook}\label{sec:conclusion}

This paper introduced MMGKS-OF, a joint image reconstruction and motion estimation technique integrated within the MMGKS framework. Numerical experiments across diverse dynamic computerized tomography applications have demonstrated the effectiveness of MMGKS-OF in reconstructing high-quality images, particularly in challenging scenarios involving limited data or significant motion.

Currently, MMGKS-OF relies on the optical flow assumption, limiting its applicability to 2D scenarios with small time differences and adherence to brightness constancy. We therefore aim to generalize our motion models beyond optical flow to encompass nonlinear deformations, as explored in works such as \cite{chen2019new, chen2021spatiotemporal}, and to develop more sophisticated models capable of handling larger time differences and 3D imaging. We also mention that small changes in time are crucial for the success of this method as it relies on the optical flow assumption. Large changes may require more complicated and complex models.

Solving the optical flow equations to estimate the parameters of the state model is computationally expensive or even prohibitive for large-scale problems. Future research will investigate techniques to reduce this cost. These include learning a reduced-order model for $\mathcal{M}$ (instead of the full model), employing machine learning approaches such as the Sparse Identification of Nonlinear Dynamics (SINDy) algorithm, and exploring model learning within the framework of linear transport maps.

The problem that we solve in this manuscript, initially introduced in \eqref{eq: dynamicEq}--\eqref{eq: measurementEq}, has strong similarities to state-space models considered in data assimilation (DA) \cite{solonen2016dimension,spantini2018inference}, but the setting we consider for dynamic tomography rely on indirect observations, described by the action of an ill-posed forward operator $\bH$, rather than direct pointwise evaluations of the state. We foresee extending this work in the Kalman filtering and smoothing setting, where promising results for smoothing (Tikhonov-type) priors are shown in \cite{hakkarainen2019undersampled} for the undersampled dynamic tomography problem.

Finally, a more rigorous theoretical analysis of the proposed algorithm's convergence properties is warranted to provide deeper insights into its behavior and potential limitations. Addressing these aspects will further enhance the performance and applicability of MMGKS-OF.

\section{Acknowledgements}

M.P. acknowledges partial support from NSF Award No. 2202846 and support from NSF DMS No. 2410699. T.O. and M.K. acknowledge support from NSF DMS No. 2410698. Any opinions, findings, conclusions, or recommendations expressed in this material are those of the authors and do not necessarily reflect the views of the National Science Foundation.
M.P. would like to further acknowledge funding provided by the Deutsche
Forschungsgemeinschaft (DFG) - Project-ID 318763901 - SFB1294 to visit the University of Potsdam where initial work on this manuscript was undertaken. M.P. would like to thank Felix Lucka for initial discussions and pointing out to resources regarding optical flow estimation for dynamic inverse problems.
\appendix

\section{Investigating Parameter Sensitivity}
Accurately selecting regularization parameters is a well-known challenge in motion estimation and image reconstruction, particularly for algorithms employing alternating minimization, as discussed in \ref{sec:mm-gks-of}. This choice often dictates the quality of the reconstruction. To demonstrate the sensitivity of a related method to its regularization parameters, we conducted an experiment using the MATLAB code and pinball dataset described in Test 4 above and we used Matlab code publicly available on the github repository \newline \url{https://github.com/HendrikMuenster/JointMotionEstimationAndImageReconstruction.git} associated to \cite{Dir15}. Similar to our proposed method, their approach employs regularization through the utilization a primal-dual framework requiring manual parameter selection. While capable of producing satisfactory results, this manual tuning presents a significant drawback necessitating the method to be run several times to estimate the desired regularization parameter. In contrast, our MMGKS method inherently performs automatic regularization parameter selection. Manual parameter tuning becomes especially impractical in large-scale problems due to increased complexity and computational cost.

\begin{figure}[htbp]
    \centering
    \begin{subfigure}{0.45\textwidth}
        \centering
    \includegraphics[width=\linewidth]{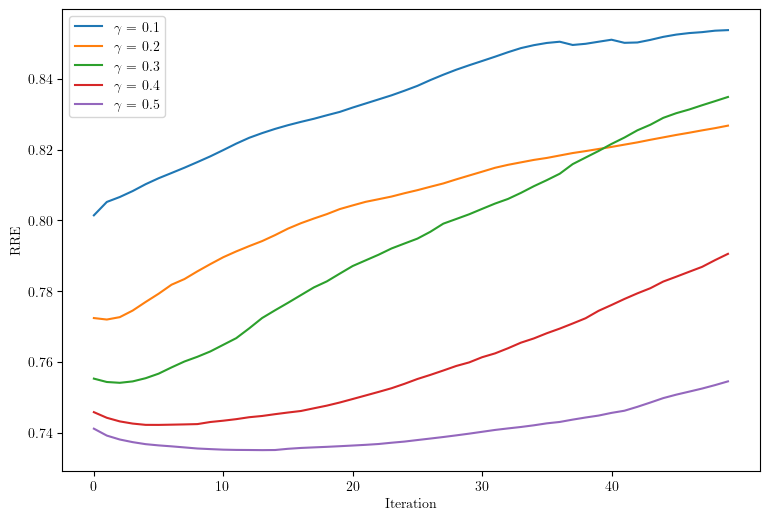}
    \caption{RRE convergence for different values of $\gamma$.}
    \label{fig:vary_gamma_rre_convergence}
    \end{subfigure}%
    \hspace{0.01\textwidth} 
    \begin{subfigure}{0.45\textwidth}
        \centering
    \includegraphics[width=\linewidth]{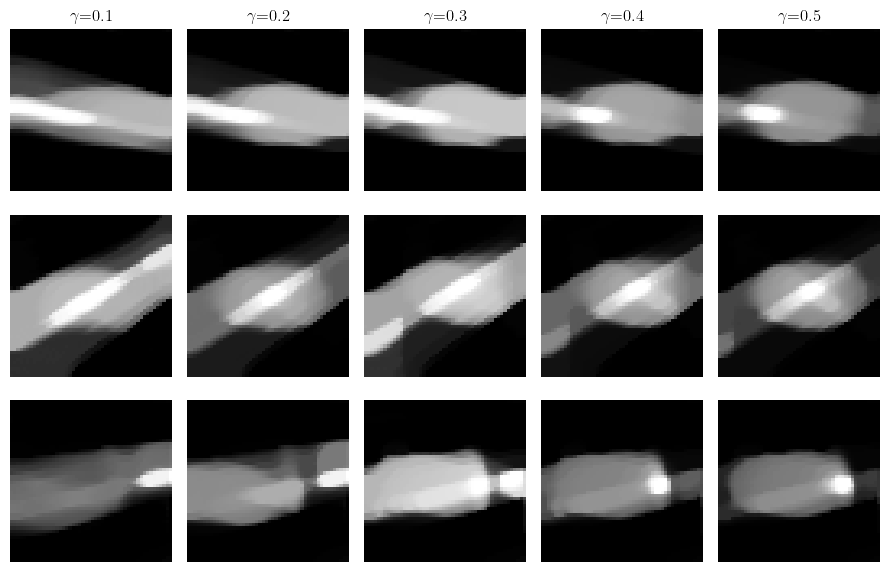}
    \caption{Reconstructed Images at $t = 2, 14, 28$ for different values of $\gamma$.}
    \label{fig:vary_gamma_rec_images}
    \end{subfigure}%
    \caption{Sensitive Dependence on Parameter $\gamma$}
    \end{figure}

Specifically, we investigated the influence of the parameter $\gamma$ on their algorithm's performance. Using values of $p$ and $q$ consistent with our own setting, we executed their MATLAB code, varying one of the regularization parameters $\gamma$ and observing the convergence of the Relative Reconstruction Error (RRE). The resulting RRE convergence curves are shown in Figure \ref{fig:vary_gamma_rre_convergence}, and the corresponding reconstructed images for different $\gamma$ values are presented in Figure \ref{fig:vary_gamma_rec_images}.

\section{MM-GKS Algorithm} 
\label{sec:mmgks_alg}
Here we present the MM-GKS algorithm that was introduced in \cite{lanza2015generalized}.
\begin{algorithm}[H]
\caption{MM-GKS \cite{lanza2015generalized}}
\label{alg:mm-gks}
\begin{algorithmic}[1]
\Require: $\bA $, $\bL $,$\bb$, ${\bu }^{(0)}$, $\epsilon$
\Ensure: An approximate solution ${\bu }^i$
\STATE \textbf{function} \textsc{MM-GKS}(${\bA }$, $\bL $, $\bb$)
\STATE Generate the initial subspace basis $\bV_\ell \in \R^{N \times \ell}$ such that $\bV_\ell^T\bV_\ell = \bI $
\FOR {$k = 0, 1, 2, \ldots$ until convergence}
\STATE $\by^{(k)} =  \bL  \bu ^{(k)}$
\STATE $\bw_\epsilon^{(k)} = ((\by^{(k)})^2 + \epsilon^2)^{1/2} $
\STATE $\bP_\epsilon^{(k)}= (diag(\bw_\epsilon^{(k)}))^{\tfrac{1}{2}}$
\STATE Obtain $\bA \bV_{\ell+k}$ and $\bP_\epsilon^{(k)}\bL \bV_{\ell+k}$
\STATE $\bA  \bV_{\ell+k} = \bQ_\bA \bR_\bA $ and $\bP_\epsilon^{(k)}\bL \bV_{\ell+k} = \bQ_\bL \bR_\bL $ \Comment{Compute/Update the QR decomposition } 
\STATE Select $\lambda^{(k)}$(e.g. GCV).
\STATE Solve the least squares problem with $\lambda^{(k)}$ to obtain $\bz^{(k+1)}$
\STATE $\bu ^{(k+1)} = \bV_{\ell+k}\bz^{(k+1)}$ 
\STATE $\br^{(k+1)} = \bA ^T(\bA \bV_{\ell+k}\bz^{(k+1)} - \bb) +\lambda^{(k)}\bL ^T \bP_\epsilon^{(k)}\bL \bV_{\ell+k}\bz^{(k+1)}$ 
\STATE $\br^{(k+1)}= \br^{(k+1)} - \bV_{\ell+k}^T\br^{(k+1)}$
\STATE $\bv = \tfrac{\br^{(k+1)}}{\|\br^{(k+1)\|_2}}$
\STATE $\bV_{\ell+k+1} = [\bV_{\ell+k}, \bv]$ \Comment{Enlarge the solution subspace} 
\ENDFOR
\STATE \textbf{end function}
\end{algorithmic}
\end{algorithm}

\section{Computing the Gradients in the OF Algorithm} 
\label{sec:of_gradients}
In the optical flow equations, we need the gradient of the image with respect to $x$, $y$ and $t$. We explain in this section how we compute this.
In other to minimize errors, we find the average around a pixel. For a sequence of $n_x \times n_y$ images $\{\bu(x,y,t)\}$, we obtain the partial derivative w.r.t time at $(x,y,t)$ by finding the average forward difference for $(2 \Delta x + 1)(2 \Delta y + 1)$ pixels around the pixel at location $(x^i,y^i)$.

\scriptsize
\begin{align}
    \frac{\partial \bu(x^i,y^i,t)}{\partial t} & = \frac{1}{(2 \Delta x + 1)(2 \Delta y + 1)} \sum_{j= -\Delta x}^{\Delta x}\sum_{k= -\Delta y}^{\Delta y} \frac{\bu(x^i+j,y^i+k,t+\Delta t) - \bu(x^i+j,y^i+k,t-\Delta t)}{2\Delta t} \\
    \frac{\partial \bu(x^i,y^i,t)}{\partial x} & = \frac{1}{2 \Delta x + 1} \sum_{j= -\Delta x}^{\Delta x} \frac{\bu(x^i+\Delta x,y^i+j,t) - \bu(x^i - \Delta x,y^i+j,t)}{2\Delta x} \\
    \frac{\partial \bu(x^i,y^i,t)}{\partial y} & = \frac{1}{2 \Delta x + 1} \sum_{j= -\Delta y}^{\Delta y} \frac{\bu(x^i+j,y^i+\Delta y ,t) - \bu(x^i+j ,y^i- \Delta y,t)}{2 \Delta y}
\end{align}

\normalsize
% \textcolor{red}{Add a comment how robust is the time direction with respect to the finite differences.}
For pixels at the edge of the image, the image is flipped to complete the box surrounding the pixels. 
For now, we set $\Delta r := \Delta x = \Delta y$ and adjust based on the magnitude of the movements (we found it effective to increase $\Delta r$ as the movement magnitude increases). 

For the velocity's gradient operator in the regularized problem, we use the finite difference discretization of the first derivative operator.

\bibliographystyle{siamplain}
\bibliography{Arxiv}
\end{document}